\documentclass[12pt]{article}
\usepackage{amsmath}
\usepackage{amssymb}

\usepackage{graphicx}

\newcommand{\qed}{\hfill$\Box$\par\medskip\par\relax}

\newcommand{\1}[1]{{\mathbf 1}{\{#1\}}}

\newcommand{\vr}{\varrho}

\newcommand{\nuo}{\nu^\omega}

\newcommand{\eps}{\varepsilon}
\newcommand{\Z}{{\mathbb Z}}

\newcommand{\Sph}{{\mathbb S}}

\newcommand{\R}{{\mathbb R}}
\newcommand{\RR}{{\mathcal R}_\omega{}}

\newcommand{\BB}{{\mathcal B}}

\newcommand{\EE}{{\mathcal E}}

\let\phi=\varphi

\newcommand{\8}{{\infty}}

\newcommand{\IQ}{{\mathbb Q}}

\newcommand{\IP}{{\mathbb P}}
\newcommand{\IE}{{\mathbb E}}
\newcommand{\Po}{{\mathtt P}_{\!\omega}}
\newcommand{\Eo}{{\mathtt E}_\omega}
\newcommand{\Eozz}{{\mathtt E}_{\omega,\zeta'}}
\newcommand{\Pozz}{{\mathtt P}_{\omega,\zeta'}}

\newcommand{\Pozx}{{\mathtt P}_{\!\omega,\zeta}}
\newcommand{\Eozx}{{\mathtt E}_{\omega,\zeta}}

\newcommand{\UU}{{\mathcal U}}

\newcommand{\nn}{{\mathbf n}_\omega{}}

\newcommand{\e}{\mathbf{e}}

\newcommand{\essinf}{\mathop{\mathrm{ess\,inf}}}

\newcommand{\near}
{\stackrel{{\scriptstyle \omega}}{\leftrightarrow}}

\newcommand{\intl}{\int\limits}

\allowdisplaybreaks

\newtheorem{theo}{Theorem}[section]
\newtheorem{lmm}[theo]{Lemma}
\newtheorem{df}[theo]{Definition}
\newtheorem{prop}[theo]{Proposition}

\newtheorem{rem}[theo]{Remark}

\title{Ballistic regime for random walks in random environment 
with unbounded jumps and 
        Knudsen billiards}
\author{Francis~Comets$^{1,}$\thanks{Partially supported by
CNRS (UMR 7599 ``Probabilit{\'e}s et Mod{\`e}les Al{\'e}atoires'')}
\and Serguei~Popov$^{2,}$\thanks{Partially supported by FAPESP
(2009/52379--8, 2009/08665--6), PRONEX ``Probabilidade e Modelagem
Estoc\'astica'', and CNPq (300886/2008--0, 472431/2009-9)}}

\begin{document}

\maketitle

{\footnotesize
\noindent $^{~1}$Universit{\'e} Paris Diderot (Paris 7), 
UFR de Math{\'e}matiques,
case 7012, Site Chevaleret,
75205 Paris Cedex 13, France\\
\noindent e-mail: \texttt{comets@math.jussieu.fr},
\noindent url: \texttt{http:/\!/www.proba.jussieu.fr/$\sim$comets}

\smallskip
\noindent $^{~2}$Department of Statistics,
Institute of Mathematics, Statistics and Scientific Computation,
University of Campinas--UNICAMP,
rua S\'ergio Buarque de Holanda 651, 13083--859, 
Campinas, SP, Brazil\\
\noindent e-mail: \texttt{popov@ime.unicamp.br}, 
\noindent url: \texttt{http:/\!/www.ime.unicamp.br/$\sim$popov}

}

\begin{abstract}
We consider a random walk in a stationary ergodic environment 
in~$\Z$, with unbounded jumps.  
In addition to uniform ellipticity and a bound on the tails of the 
possible jumps, we assume a condition of strong transience to the
right which implies that there are no ``traps''.  
We prove the
law of large numbers with positive speed, as well as the ergodicity
of the environment seen from the particle. Then, we consider
Knudsen stochastic billiard with a drift in a 
random tube in~$\R^d$, $d\geq 3$, which serves as environment.
The tube is infinite in the first direction, 
and is a stationary and ergodic process indexed by the first 
coordinate.
A particle is moving in straight line inside the tube, 
and has random bounces upon hitting the boundary, according to 
the following modification of the cosine reflection law: 
the jumps in the positive direction are always
accepted while the jumps in the negative direction may be rejected.
Using the results for the random walk 
in random environment together with an appropriate coupling,
we deduce the law of large numbers for the stochastic billiard 
with a drift.
\\[.3cm]\textbf{R\'esum\'e:} Nous consid\'erons une marche al\'eatoire dans un milieu stationnaire ergodique sur $\Z$, 
avec des sauts non born\'es. En plus de l'uniforme ellipticit\'e et d'une borne uniforme sur la queue de la loi des sauts,
nous supposons une condition de transience forte qui garantit l'absence de "pi\`eges". Nous montrons la loi des grands nombres avec vitesse strictement positive, ainsi que l'ergodicit\'e de l'environnement vu de la particule. Par ailleurs, nous
\'etudions aussi le billard stochastique de Knudsen avec d\'erive dans un tube al\'eatoire dans $\R^d$, $d \geq 3$, 
qui constitue l'environnement.  Le tube est infini dans la premi\`ere direction, et, vu comme un processus ind\'ex\'e par la premi\`ere coordonn\'ee, il est suppos\'e stationnaire ergodique.
 Une particule se d\'eplace en ligne droite \`a l'int\'erieur du tube, avec des rebonds al\'eatoires sur le bord, 
selon la modification suivante de la loi de reflexion 
en cosinus: les sauts dans la direction positive sont toujours accept\'es, tandis que ceux dans l'autre direction peuvent \^etre rejet\'es. En utilisant les r\'esultats pour la   marche al\'eatoire en milieu al\'eatoire et un couplage appropri\'e, 
nous obtenons la   loi des grands nombres pour le  billard stochastique avec d\'erive.
 \\[.3cm]\textbf{Keywords:} cosine law, stochastic billiard,
Knudsen
random walk,
random medium, random walk in random environment, unbounded jumps,
%non-i.i.d.\ 
stationary ergodic environment, regenerative structure,  
point of view of the particle
\\[.3cm]\textbf{AMS 2000 subject classifications:}
 60K37. Secondary: 37D50, 60J25
\end{abstract}

\section{Introduction}
\label{s_intro}

Stochastic billiards deal with the motion of a particle 
inside a connected domain in the Euclidean space, travelling in
straight lines inside the domain and subject to random bouncing 
when hitting the boundary. They are motivated by problems 
of transport and diffusion inside nanotubes, where the complex
microscopic structure of the tube boundary allows for 
a stochastic description of the collisions: they can be viewed as 
limits of deterministic billiards on tables with rough boundary
as the ratio of macro to micro scales diverges~\cite{Feres07}.
See also~\cite{CPSV1} and~\cite{CPSV3} for a detailed perspective
 from physics and chemistry, and~\cite{Feres07, CPSV1} 
%%for basic results in the case of bounded domains in the first
% reference.
for basic results.
A natural reflection law is when the outgoing direction has a 
density proportional to the cosine of its angle with the 
inner normal vector,
independently from the past. This model, originally introduced by
Martin 
Knudsen, is called the Knudsen stochastic billiard. It
has two important features: the uniform measure on the phase space 
is invariant for the dynamics, and moreover it is reversible. To
understand large time behavior of Knudsen billiards, one needs to
consider infinite domains. Recurrence and transience is studied 
in~\cite{MVW}, for billiards in a planar tube extending to infinity
in the horizontal direction, under assumptions of regularity and
growth on the tube. For the physically relevant case of an
infinite tube   which is irregular but has some homogeneity
properties at 
large scale, the description of the tube as a random environment 
has been introduced in~\cite{CPSV2}: the domain is the realization
of a stationary ergodic process indexed by the horizontal
coordinate.
Diffusivity of the particle is studied in this paper in dimension
$d=1+(d-1) \geq 2$: generically, when the tube does not have 
arbitrarily long cavities, the billiard is diffusive in
dimension $d \geq 3$, and also for $d=2$ when the
billiard has ``finite horizon''. 
Reversibility allows to find the limit of 
the environment seen from the particle, and to use the appropriate
techniques  which have been extended to Random Walks in Random
Environments (RWRE). We briefly mention~\cite{CPSV3} for
the non-equilibrium dynamics aspects of the billiard 
and some features as a microscopic model for diffusion.

In the model of~\cite{CPSV2} the large-scale picture of the motion
of the particle is purely diffusive; in particular, the limiting
velocity of the particle equals zero.
In this paper, we consider a stochastic billiard in a random tube
as in~\cite{CPSV2} 
% but with a constant drift in the reflection law, 
traversed by a flow with constant current to the right;
our goal is to prove the law of large numbers (with positive 
limiting velocity). 
This current is modelled in the following way: the jumps in the
positive direction are always accepted,
but the jumps in the negative direction are accepted with
probability~$e^{-\lambda u}$, where~$u$ is the horizontal size of
the attempted jump. This method of giving a drift to the particle
has the following advantage: the reversibility of the 
stochastic billiard is preserved (although, of course, the
reversible measure is no longer the same), which simplifies
considerably
the analysis of the model.
%  Thus,
% we consider Knudsen billiard in a random tube as in~\cite{CPSV2},
% but with a constant drift in the reflection law. 
In view of the above, the large scale picture is expected 
to be similar to the one-dimensional RWRE with a drift
 when the environment is given by a ``resistor network'' with a
similar acceptance/rejection mechanism; 
the environment is not i.i.d.\ but stationary ergodic, the jumps
are not nearest neighbor but \textit{unbounded}.

We review known results on the law of large numbers for transient 
RWREs on~$\Z$. For nearest neighbor jumps, the sub-ballistic and
ballistic regimes -- meaning that the speed is zero, 
resp.\ non zero -- are fully understood 
(e.g., \cite{Zeit} and Section~1 in~\cite{Sz-tr}) 
with the explicit formula of Solomon for the speed~\cite{Sol}
in the case of i.i.d. environment; the extension to stationary 
ergodic environment is given in~\cite{Alili}. When
the  jumps are bounded but not nearest neighbor, 
Key~\cite{Key} shows that transience of the walk to the
right amounts to positivity of some middle Lyapunov
exponent of a product of random matrices. The regime
where the law of large numbers holds with a positive speed is
characterized  in Br\'emont's~\cite{Bre} 
by the positivity of this exponent and existence of 
an invariant law for the environment absolutely continuous to the 
static law,
%(Theorem 1.10). 
but no explicit formula is anymore available. 
Goldsheid~\cite{Gold} gives
sufficient conditions 
(which are also necessary in the case of i.i.d.\
environment), and also for the quenched central limit theorem. 
For completeness, we mention a result of 
Bolthausen and Goldsheid~\cite{BoGo}
for recurrent RWREs with bounded jumps: if the quenched drift is
not a.s.\ zero, the typical displacement at time~$n$ 
is of order $\ln^2 n$,
i.e., 
the RWRE has a similar lingering behavior as Sina\"{\i}'s walk.

The case of unbounded jumps has been very seldom considered; 
in fact, we can only mention that Andjel~\cite{A} proves 
a 0-1 law when the jumps
have uniform exponential tails.

In this paper, we prove the law of large numbers with a positive 
speed for RWRE on~$\Z$ with unbounded jumps, 
under the following assumptions: stationary ergodic environment, 
(E) uniform ellipticity; 
(C) uniform (and integrable, but not necessarily exponential) tails
for the jumps; (D) strong uniform transience to the right.
We do not assume reversibility of the RWRE. The
strategy is to consider an auxiliary RWRE with truncated jumps, 
to prove the existence of limits for the speed and the environment 
seen from the walker, then let the truncation parameter 
tend to infinity, and find a limit point for the environment 
measure. We mention also that assumption~(D) precludes the 
existence of arbitrarily long traps -- i.e., pieces of the
environment where the random walk can spend an unusually large time
--, and it is rather strong. We emphasize that we do not assume any
mixing -- hence, no independence -- on the environment.
As we see below, this set of conditions is adapted to  our purpose.
In our opinion,
it is a challenging problem to find weaker conditions that still
permit to obtain the law of large numbers for RWREs with unbounded
jumps with only polynomial tails. In particular, it would be especially
interesting to substitute the current condition~(D) by a weaker one; however,
at the moment we do not have any concrete results and/or plausible conjectures
which go in that direction.

To apply this result to the billiard in random tube, we need a
discretization 
procedure to compare the billiard to a random walk. 
 This can be performed by coupling the billiard 
with an independent coin tossing; the integer part
of the horizontal coordinate of hitting points on the boundary,
sampled at success times of the coin tossing, is an embedded RWRE 
in some environment determined by the random tube. 
We check condition~(D) for the RWRE by making use of the
reversibility of the billiard and spectral estimates 
(as in~\cite{Sh} for a reversible RWRE on~$\Z^d$).
This coupling allows to transfer results from the  
RWRE -- a simplified model -- to the stochastic billiard
-- a much more involved one --.
Under fairly reasonable assumptions on the random tube, 
we obtain for the billiard the law of large numbers with positive
speed.

The paper is organized as follows: 
we define the two models and state the results in the next section.
Section~\ref{s_proofs_rwre} contains the proofs for RWRE, and
Section~\ref{s_proofs_billiard} those for the
stochastic billiard, including the construction of 
the coupling with the RWRE.

\section{Formal definitions and results}
\label{s_results}
Now, we formally define the random billiard with drift in a random 
tube and the one-dimensional random walk in 
stationary ergodic random environment with unbounded jumps.

Already at this point we warn the reader that the (continuous)
random environment for the billiard processes and the (discrete)
random environment for the random walk are denoted by the same
letter~$\omega$. Hopefully, this creates no confusion since at all
times we tried to make it clear which model is under consideration.
Also, we use the same notation~$\IP$ (the law of the random
environment) 
% and $\IQ$ (the law of, possibly truncated,
% environment-seen-from-the-particle process) 
for both models.

\subsection{Random billiards with drift}
\label{s_results_billiards}
We define the model of random billiard in a random tube, basically
keeping the notations of~\cite{CPSV2}.

In this paper, $\R^{d-1}$ will always stand for the linear
subspace of $\R^d$ which is perpendicular to the first coordinate
vector~$\e$,
we use the notation~$\|\cdot\|$ for the Euclidean norm in~$\R^d$.
% or~$\R^{d-1}$. 
%For $k\in \{d-1,d\}$
Let $\BB(x,\eps)=\{y\in\R^d:\|x-y\|<\eps\}$ be the open
$\eps$-neighborhood of $x\in\R^d$. Define
 $\Sph^{d-1}=\{y\in\R^d:\|y\|=1\}$ to be the unit sphere in~$\R^d$.
% and let $\Sph^{d-2}=\Sph^{d-1}\cap \R^{d-1}$ be the unit sphere
% in~$\R^{d-1}$.
We write $|A|$ for the $d$-dimensional Lebesgue measure in case
$A\subset \R^d$,
and $(d-1)$-dimensional Hausdorff measure in case 
$A\subset \Sph^{d-1}$. Let
\[
\Sph_h = \{w\in\Sph^{d-1}: h\cdot w > 0\}
\]
be the half-sphere looking in the direction~$h$.
 For $x\in\R^d$, it will frequently be convenient to write
$x=(\alpha,u)$, being~$\alpha$ the first coordinate of~$x$ and
$u\in\R^{d-1}$; then, $\alpha=x\cdot \e$, and we write $u=\UU x$,
where~$\UU$ is the projector on~$\R^{d-1}$. 
Fix some positive constant~${\widehat M}$, and define
\begin{equation}
\label{def_Lambda}
 \Lambda = \{u\in \R^{d-1} : \|u\| \leq {\widehat M}\}.
\end{equation}

%%Let~$A$ be an open connected domain in~$\R^d$. 
We denote 
by~$\partial A$ the boundary of~$A \subset \R^d$,
by~$\bar A = A\cup \partial A$ the closure of~$A$
and by $A^\circ$ the interior of~$A$ (i.e., the largest open
set contained in~$A$).

\begin{df}
\label{def_Lipschitz}
 Let $k\in \{d-1,d\}$, and 
$A$ a subset of~$\R^k$. 
%% suppose that~$A$ is  an open domain in~$\R^k$. 
We say that~$\partial A$ is
$({\hat\eps},{\hat L})$-Lipschitz, if for any~$x\in\partial A$ 
there exist
an affine isometry ${\mathfrak I}_x : \R^k\to\R^k$ and a function
$f_x:\R^{k-1}\to\R$ such that
\begin{itemize}
\item $f_x$ satisfies Lipschitz condition with 
constant~${\hat L}$, 
i.e., $|f_x(z)-f_x(z')| \leq {\hat L}\|z-z'\|$ for all $z,z'$;
\item ${\mathfrak I}_x x = 0$, $f_x(0)=0$, and
\[
%  {\mathfrak I}_x(A\cap\BB(x,{\hat\eps})) 
%             = \{z\in\BB(0,{\hat\eps}): 
%                       z^{(k)} > f_x(z^{(1)},\ldots,z^{(k-1)})\}.
  {\mathfrak I}_x\big(A^\circ \cap\BB(x,{\hat\eps})\big) 
             = \big\{z\in\BB(0,{\hat\eps}): 
                       z^{(k)} > f_x(z^{(1)},\ldots,z^{(k-1)})\big\}.  
\]
\end{itemize}
% In the degenerate case $k=1$ we adopt the convention 
% that~$\partial A$ is
% $({\hat\eps},{\hat L})$-Lipschitz for any positive
%  ${\hat\eps},{\hat L}$.
\end{df}

Now, fix ${\widehat M}$, and define~$\EE$
to be the set of all open domains~$A$ such that $A\subset\Lambda$
and~$\partial A$ is $({\hat\eps},{\hat L})$-Lipschitz for some
$({\hat\eps},{\hat L})$ (which may depend on~$A$).
We turn~$\EE$ into a metric space by defining the distance 
between~$A$ and~$B$ to be equal to
$|(A\setminus B)\cup(B\setminus A)|$.
Let $\Omega$
%=\cCC_\EE(\R)$ 
be the space of all 
%continuous
c\`adl\`ag
functions $\R\to\EE$, let~${\mathcal A}$ be the sigma-algebra
generated by the cylinder sets with respect to the Borel
sigma-algebra
on~$\EE$, and let~$\IP$ be a probability measure on 
$(\Omega,{\mathcal A})$. 
This defines
a $\EE$-valued process $\omega=(\omega_\alpha, \alpha\in\R)$.
Write~$\theta_\alpha$ for the spatial shift: $\theta_\alpha
\omega_{\cdot} = \omega_{\cdot+\alpha}$. We suppose that the
process~$\omega$
is stationary and ergodic with respect to the family of shifts
$(\theta_\alpha,\alpha\in\R)$.
With a slight abuse of notation, we denote also by
\[
 \omega = \{(\alpha,u)\in\R^d : u\in \omega_\alpha\}
\]
the random domain (``tube'') where the billiard lives. 
Intuitively, $\omega_\alpha$ is the ``slice'' obtained 
by crossing~$\omega$ with the hyperplane 
$\{\alpha\}\times \R^{d-1}$. 
% One can check that the domain~$\omega$ is an open subset of~$\R^d$, and we will assume that it is connected.

We will assume that the domain~$\omega$ is connected. A trivial sufficient condition  is that $\omega_\alpha$ is connected for all $\alpha$; a typical example is when $\partial \omega$ is generated by rotating around the horizontal axis the graph of a one-dimensional (stationary ergodic) process with values in $[1, \widehat M]$. In this paper, we will work under the
more general Condition~P below, which implies that $\omega$ is arc-connected.

We also assume the following

\medskip
\noindent
\textbf{Condition~L.}
There exist ${\tilde\eps},{\tilde L}$ such that $\partial\omega$
is $({\tilde\eps},{\tilde L})$-Lipschitz (in the sense of
Definition~\ref{def_Lipschitz}) $\IP$-a.s.

\medskip

Denote by~$\nuo$ the $(d-1)$-dimensional Hausdorff measure
on $\partial\omega$; from Condition~L one obtains that~$\nuo$
is locally finite.
We keep the usual notation $dx, dv, dh, \ldots$ for the
$(d-1)$-dimensional Lebesgue measure
on~$\Lambda$ (usually restricted to $\omega_\alpha$ for
some~$\alpha$) or the
surface  measure
on~$\Sph^{d-1}$.

% \begin{figure}
% \centering
% \includegraphics[width=\textwidth]{tube}
% \caption{Notations for the random tube.
% Note that the sections may be
% disconnected}
% \label{f_tube}
% \end{figure}

Define the set of regular points 
\[
 \RR = \{x\in \partial\omega : \partial\omega \text{ is
continuously differentiable in }x \}.
\]
For all $x=(\alpha,u)\in\RR$, let us 
define also the normal vector $\nn(x)=\nn(\alpha,u)\in\Sph^{d-1}$ 
pointing inside the domain~$\omega$.

We suppose that the following condition holds:

\medskip
\noindent
\textbf{Condition~R.} 
We have $\nuo(\partial\omega\setminus\RR)=0$, $\IP$-a.s.

\medskip

We say that $y\in\bar\omega$ is \emph{seen from} 
$x\in\bar\omega$ if there exists $h\in\Sph^{d-1}$ and~$t_0>0$ such
that $x+th\in\omega$ for all $t\in (0,t_0)$ and $x+t_0 h = y$. 
Clearly, if~$y$ is seen from~$x$ then~$x$ is seen from~$y$,
and we write ``$x \near y$'' when this occurs. 

\begin{figure}
 \centering
\includegraphics{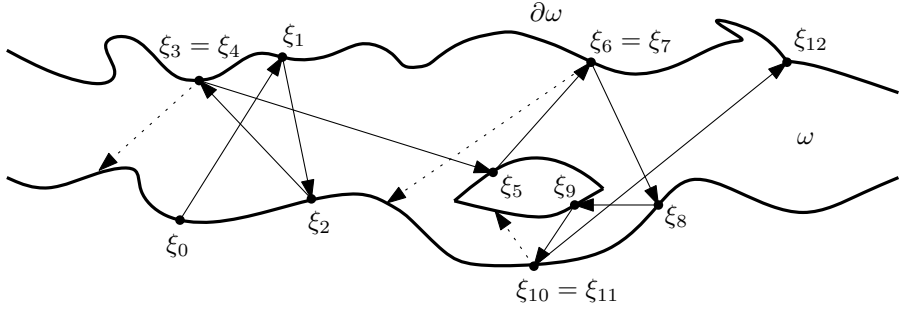}
\caption{Knudsen random walk with drift}
\label{f_bil_drift}
\end{figure}

One of the main objects of study in the paper~\cite{CPSV2} is the 
Knudsen random walk (KRW) which is a
discrete time Markov process on~$\partial\omega$,
defined through its transition density~$K$ with respect to the surface measure $\nu^\omega$:
for $x,y\in\partial\omega$
\begin{equation}
\label{def_K}
 K(x,y) = \gamma_d
\frac{\big((y-x)\cdot\nn(x)\big)\big((x-y)\cdot\nn(y)\big)}
{\|x-y\|^{d+1}} \1{x,y\in\RR, x \near y},
\end{equation}
where $\gamma_d = \big(\int_{\Sph_{\e}} h\cdot \e\, dh\big)^{-1}$
is the normalizing constant. We also refer to the Knudsen 
random walk as the random walk with cosine
reflection law, since it can be easily seen from~\eqref{def_K} that
the density of the outgoing direction is proportional to the cosine
of the angle between this direction and the normal vector (see,
 e.g., formula~(4) in~\cite{CPSV1}). In this paper, however, we
shall consider the walk which ``prefers'' the positive direction: a
jump in the direction~$\e$ is always accepted, 
but if the walk attempts to
jump in the negative direction~$(-\e)$, it is accepted with
probability
$e^{-\lambda u}$, where~$u$ is the horizontal size of the attempted
jump and~$\lambda>0$ is a given parameter. Formally, define
\begin{equation}
\label{def_hat_K}
 {\hat K}(x,y) = \begin{cases}
                  K(x,y), & \text{if } (y-x)\cdot \e \geq 0,\\
  e^{\lambda (y-x)\cdot \e}K(x,y), & \text{if } (y-x)\cdot \e < 0,
                 \end{cases}
\end{equation}
and let
\begin{equation}
\label{def_Theta}
 \Theta(x) = 1-\intl_{\partial\omega} {\hat K}(x,y)\,d\nuo(y).
\end{equation}
%%Being $\Po,\Eo$ the quenched (i.e., with fixed~$\omega$)
%%probability and expectation, we define the 
%%\emph{Knudsen random walk with drift} (KRWD) in the 
%%following way:
For a fixed~$\omega$, we define the \emph{Knudsen random walk with
drift} (KRWD) $(\xi_n; n \geq 0)$ -- and denote by   $\Po,\Eo$ the
corresponding quenched probability and expectation -- as the Markov
chain on ${\partial\omega}$ starting from $\xi_0=0$ such that,
for any $x\in\RR$ and any measurable
$B\subset \partial\omega$ such that $x\notin B$,
\[
\Po[\xi_{n+1}\in B \mid \xi_n=x] = \intl_B {\hat K}(x,y)\,d\nuo(y),
\]
and
\[
 \Po[\xi_{n+1}=x \mid \xi_n=x] = \Theta(x).
\]
On Figure~\ref{f_bil_drift} one can see a typical path of the
random walk (rejected jumps are shown as dotted lines).

 As observed in~\cite{CPSV1},
$K(\cdot,\cdot)$ is symmetric (that is, $K(x,y)=K(y,x)$ for all
$x,y\in\RR$), so that the $(d-1)$-dimensional Hausdorff
measure~$\nuo$ is reversible for~$K$. Then, with
$\pi(x)=e^{\lambda(x\cdot\e)}$, the
measure~$\nuo_\lambda$ on~$\partial\omega$ given by
\begin{equation}
\label{def_nu_lambda}
 \frac{d\nuo_\lambda}{d\nuo}(x) = 
\pi(x)
%%\stackrel{\rm def}{=}
=e^{\lambda(x\cdot\e)}
\end{equation}
is such that $\pi(x){\hat K}(x,y)=\pi(y){\hat K}(y,x)$, 
showing that~$\nuo_\lambda$ is reversible for the KRWD~$\xi$. 
With some abuse of
notation, we shall sometimes write $\pi(B):=\nuo_\lambda(B)$ for
$B\subset\partial\omega$.

We need also require a last technical assumption:

\medskip
\noindent
\textbf{Condition~P.}
There exist constants~$N,\eps,\delta$ such that for $\IP$-almost
every~$\omega$, for any $x,y\in\RR$ with $|(x-y)\cdot \e|\leq 2$
there exist $B_1,\ldots,B_n\subset\partial\omega$,
$n\leq N-1$ with $\nuo(B_i)\geq \delta$ for all $i=1,\ldots,n$, 
and such that 
\begin{itemize}
 \item $K(x,z)\geq \eps$ for all $z\in B_1$,
 \item $K(y,z)\geq \eps$ for all $z\in B_n$,
 \item $K(z,z')\geq \eps$ for all $z\in B_i$, $z'\in B_{i+1}$,
$i=1,\ldots,n-1$
\end{itemize}
(if $N=1$ we only require that $K(x,y)\geq \eps$).
In other words, there exists a ``thick'' 
path of length at most~$N$
joining~$x$ and~$y$. This assumption is already used 
in~\cite{CPSV2}, it prevents the tube from splitting into separate
channels of arbitrary length, which could slow down the
homogeneization.
  
\medskip

% Now, following~\cite{CPSV1}, we define also the Knudsen 
% stochastic
% billiard with drift (abbreviated KSBD)
% $(X,V)$. First, we do that for the process starting on the
% boundary~$\partial\omega$
% from the point~$x_0\in \partial\omega$. Let
% $x_0=\xi_0,\xi_1,\xi_2,\xi_3,\ldots$
% be the trajectory of the random walk, and define
% \[
%  \tau_n = \sum_{k=1}^n \|\xi_k-\xi_{k-1}\|.
% \]
% Then, if $\tau_n<\tau_{n+1}$, for $t\in[\tau_n,\tau_{n+1})$
% we
% define 
% \begin{equation}
% \label{def_billiard}
%  X_t=\xi_n+(\xi_{n+1}-\xi_n)\frac{t-\tau_n}{\|\xi_{n+1}-\xi_n\|}.
% \end{equation}
% The quantity~$X_t$ stands for the position of the particle at
% time~$t$.
% Since~$(X_t)_{t \geq 0}$  is not a Markov process by itself, 
% we define also the c\`adl\`ag version
% of the motion direction at time~$t$:
% \begin{equation}
% \label{def_direction}
% V_t = \lim_{\eps \downarrow 0}\frac{X_{t+\eps}-X_t}{\eps}.
% \end{equation}
% Then, it holds that $V_t\in \Sph^{d-1}$, and the pair 
% $(X_t,V_t)_{t \geq 0}$
% is a Markov process. Of course, in a similar way one can define
% also
% the stochastic
% billiard starting from the interior of~$\omega$ by specifying its
% initial position~$X_0$ and initial direction~$V_0$.

We prove the existence of the speed of KRWD:
\begin{theo}
\label{t_kn_lln}
 Assume that $d\geq 3$. There exists a positive deterministic
${\hat v}$ such that for $\IP$-almost every~$\omega$
\begin{equation}
\label{lln_krw}
 \frac{\xi_n\cdot\e}{n} \to {\hat v} \qquad \text{as
 $n\to\infty$, $\Po$-a.s.}
\end{equation}
% \begin{align}
% \frac{\xi_n\cdot\e}{n} \to {\hat v} \qquad \text{as
% $n\to\infty$, $\Po$-a.s.},\label{lln_krw}\\
% \frac{X_t\cdot\e}{t} \to {\tilde v} \qquad \text{as
% $t\to\infty$, $\Po$-a.s.}\label{lln_ksb}
% \end{align}
\end{theo}
The assumption $d\geq 3$ crucially enters estimate~(\ref{<j}). 
In dimension~2,
depending on the geometry, KRWD can have large jumps, and then may
not obey the law of large numbers. Naturally, if one assumes a
strong additional condition that the size of the jumps is a.s.\
uniformly bounded (so-called finite horizon condition in the 
billiard literature), then our argument works in the case $d=2$ as
well. 
Still, we feel that Theorem~\ref{t_kn_lln} can hold in dimension~2
even without the finite horizon condition; for the proof, however,
one would need estimates on the size of the jump that are finer than the 
``uniform'' one provided by~\eqref{long_jump}. Note that in the driftless 
case we can control the average size of the jump using the explicit form of 
the stationary measure for the environment seen from the particle
(cf.\ Lemma~4.1 in~\cite{CPSV2}). Unfortunately, in the presence of the drift 
one does not obtain the stationary measure for this process in such an
explicit way, and this is the reason why the situation in dimension~2
is less clear.

\subsection{One-dimensional random walk in random environment}
\label{s_results_rwre}
Let us consider a collection of nonnegative numbers 
$\omega = (\omega_{xy}; x,y\in\Z)$, with the property $\sum_y
\omega_{xy} =1$ for all~$x$. This collection is called the
\emph{environment}, and we denote by~$\Omega$ the space of all
environments. Next, we consider a Markov chain
$(S_n, n=0,1,2,\ldots)$ with the transition probabilities
\[
 \Po^{x_0}[S_{n+1}=x+y\mid S_n=x] 
       = \omega_{xy} \text{ for all }n\geq 0, \quad
               \Po^{x_0}[S_0=x_0]=1,
\]
so that $\Po^{x_0}$ is the \emph{quenched} law of the Markov chain 
starting from~$x_0$ in the environment~$\omega$. 
Let us write~$\Po$ 
for $\Po^0$. The environment is chosen at random from the
space~$\Omega$ according to a law~$\IP$ before the
random walk starts. We denote by~$\Eo^{x_0}$ and~$\IE$ the expectations with
respect to~$\Po^{x_0}$ and~$\IP$ correspondingly. Also, we assume that
the sequence of random vectors $(\omega_{x\,\cdot}, x\in\Z)$ 
is stationary and ergodic. 

We need the following (one-sided) uniform ellipticity condition:

\medskip
\noindent
\textbf{Condition~E.} There exists $\tilde\eps$ such that 
$\IP[\omega_{01}\geq\tilde\eps]=1$.

\medskip

For any integer $\vr>1$ let us define also the ``truncated''
environment~$\omega^\vr$ by
\[
 \omega^\vr_{xy} = \begin{cases}
                      \omega_{xy}, & \text{ if } 0<|y|<\vr,\\
                        0, & \text{ if } |y|\geq \vr,\\
                       \omega_{x0}+
 \displaystyle\sum_{y:|y|\geq \vr}\omega_{xy}, & \text{ if } y=0,
                     \end{cases}
\]
and observe also that formally $\omega=\omega^\infty$. 
The truncated random walk~$S^\vr$ is then defined by
\[
 \Po^{x_0}[S^\vr_{n+1}=x+y\mid S^\vr_n=x] 
       = \omega^\vr_{xy} \text{ for all }n\geq 0, \quad
               \Po^{x_0}[S^\vr_0=x_0]=1,
\]
In words, the random walk~$S^\vr$ in the 
truncated environment~$\omega^\vr$ 
is the modification of the original random walk 
%such that, 
%if the size of the attempted jump is at least~$\vr$, 
%then this jump
%is rejected and the particle remains on the place. 
where jumps of lengths less than $\rho$ are kept, but larger jumps 
are rejected and the particle does not move.
We shall sometimes also write e.g.\
$\Po[S^{\vr_1}\in\cdot\;,S^{\vr_2}\in\cdot\;]$ meaning here the
natural coupling of two versions of the random walk with different
truncation but in the same environment. This coupling is defined in
the following way:
\begin{itemize}
\item if $S^{\vr_1}_n\neq S^{\vr_2}_n$, then $S^{\vr_1}_{n+1}$
and $S^{\vr_2}_{n+1}$ are independent given $S^{\vr_1}_{i},
S^{\vr_2}_{i}, i \leq n$;
\item if $S^{\vr_1}_n=S^{\vr_2}_n=x$, and $Y_n$ is a random variable
with $\Po[Y_n=y]=\omega_{xy}$ and independent of  $S^{\vr_1}_{i},
S^{\vr_2}_{i}, i \leq n$, then
\[
 S^{\vr_i}_{n+1}  = \begin{cases}
                     x+Y_n, & \text{if }|Y_n|<\vr_i,\\
                     x, & \text{if }|Y_n|\geq\vr_i,
                    \end{cases}
\]
for $i=1,2$.
\end{itemize}

Let us assume the following condition on the tails of the possible 
jumps of the random walks:

\medskip
\noindent
\textbf{Condition C.} There exist $\gamma_1>0$ and $\alpha>1$ 
such that for all $s\geq 1$ we have
\begin{equation}
\label{eq_tails}
 \sum_{y: |y| \geq s}\omega_{0y} \leq \gamma_1 s^{-\alpha}, 
     \qquad \text{$\IP$-a.s.}
\end{equation}

\medskip

For $I\subset\Z_+$ and $A\subset \Z$ we denote by $N^\vr_I(A)$ 
the number of visits to~$A$ of the random walk~$S^\vr$ during the
time set~$I$, i.e.,
\[
 N^\vr_I(A) = \sum_{k\in I}\1{S^\vr_k \in A}.
\]
We use the shorter notations $N^\vr_I(x)$ for $N^\vr_I(\{x\})$,
$N^\vr_k(A):=N^\vr_{[0,k]}(A)$ for the number of visits 
to~$A$ during the time interval~$[0,k]$, and
$N^\vr_k(x):=N^\vr_k(\{x\})$.

Next, we make another assumption that says, essentially, 
that the random walk is ``uniformly'' transient to the right
(i.e., there are no ``traps'').

\medskip
\noindent
\textbf{Condition D.} There is a non-increasing function~$g_1\geq 0$ with the
property
$\sum_{k=1}^\infty kg_1(k)<\infty$ and a finite $\vr_0$,
such that for all~$x\leq 0$ 
and all~$\vr\geq \vr_0$,  $\IP$-almost surely it holds that $\Eo^0
N^\vr_\infty(x)\leq g_1(|x|)$.

\medskip

With these assumptions, we can prove that the speed of the 
random walk is well-defined and positive:
\begin{theo}
\label{t_rw_lln}
 For all $\vr\in[\vr_0,\infty]$ there exists~$v_\vr > 0$ such that 
for $\IP$-a.a.~$\omega$
we have
\begin{equation}
\label{eq_rw_lln}
 \frac{S^\vr_n}{n} \to v_\vr \quad \text{ as $n\to\infty$,
$\Po$-a.s.}
\end{equation}
\end{theo}

Next, we are interested in the environment seen from the particle.
Let~$\theta_z$ be the shift to~$z$ acting on~$\omega$ in the 
following way: $(\theta_z\omega)_{xy}=\omega_{x+z,y}$. 
The process of the environment viewed
from the particle (with respect to~$S^\vr$)
 is defined by $\omega(n)=\theta_{S^\vr_n}\omega$.
\begin{theo}
\label{t_rw_env}
 For all $\vr\in [\vr_0,+\infty]$ there exists an unique invariant
measure~$\IQ^\vr$ for the process of the environment viewed from the
particle with~$\IQ^\vr \ll \IP$. Then, we have
\begin{equation}
\label{eq_rw_speed} 
 v_\vr = \int_{\Omega}\Eo^0 S^\vr_1 \, d\IQ^\vr.
\end{equation}
Moreover, for all $\vr\in [\vr_0,+\infty]$ the measure~$\IQ^\vr$
is ergodic and~$\IQ^\vr$ weakly converges
to~$\IQ^\infty$ as $\vr \to \infty$. Finally, it holds that
$v_\vr \to v_\infty$ as $\vr \to \infty$.
\end{theo}

\begin{rem}
\label{r_explicit}
In the case $\vr<\infty$ the invariant measure~$\IQ^\vr$ 
is given by the 
formula~\eqref{eq_Q} in Section~\ref{s_proofs_rwre}.
\end{rem}

\section{Proofs for RWRE}
\label{s_proofs_rwre}
Denote by $T^\vr_z = \min\{k\geq 0 : S^\vr_k\geq z\}$.
We use the simplified notation $T^\vr:=T^\vr_0$.
Let
\[
 r_x^\vr(z) = \Po^x[S^\vr_{T^\vr_z}=z]
\]
be the probability that, at moment~$T^\vr_z$,
 the (truncated) random walk is located \emph{exactly} at~$z$.
We also use the shorter notation $r_x^\vr:=r_x^\vr(0)$. Of course, 
the quantity $r_x^\vr(z)$ depends also on~$\omega$, but, for the 
sake of simplicity, we keep it this way.

% \textbf{(here comment on the strategy of the proof? Example with
% exponential tails?)}

The key fact needed in the course of the proof of our results is 
the following lemma:
\begin{lmm}
\label{l_unif_eps}
Assume Conditions~E, C, D. Then, there exists $\eps_1>0$ 
such that, $\IP$-a.s., 
 \[
  r_x^\vr\geq 2\eps_1
 \]
for all $x\leq 0$ and for
all~$\vr\in [\vr_0, \infty]$.
\end{lmm}

\noindent
\textit{Proof.}
Let us denote $\Z_-^*:=\Z\cap (-\infty,-1]$, 
 $\Z_-:=\Z\cap (-\infty,0]$.
%% there exists~$C_1$ such that for all~$x < 0$
%%\begin{equation}
%%\label{expected_hit_negative}
%% \Eo^{x}N^\vr_\infty(\Z_-) \leq C_1 |x|
%%\end{equation}
For~$\vr \geq \vr_0$ and $x \leq -1$,
\begin{align}
 \Eo^{x}N^\vr_\infty(\Z_-) 
& = \sum_{z \leq x}    \Eo^{x}N^\vr_\infty(z) + \sum_{x <z\leq 0} 
   \Eo^{x}N^\vr_\infty(z) \nonumber\\
& = \sum_{z \leq x}    \Eo^{x}N^\vr_\infty(z) + \sum_{x <z\leq 0} 
   \Eo^{z}N^\vr_\infty(z) \qquad {\rm (Markov\ property)}\nonumber\\
& \leq \sum_k g_1(k) + |x| g_1(0)  \qquad \qquad \qquad \quad {\rm (Condition~D)} \nonumber\\
\label{expected_hit_negative}
 & \leq C_1 |x|
\end{align}
for $\IP$-almost all~$\omega$, 
with some finite constant $C_1$. 
Since 
\[
 \{T^\vr>k\}= \{N^\vr_k(\Z_-^*)\geq k+1\} \subset \{N^\vr_k(\Z_-)\geq k\}
\]
for such an $x$, 
using Chebyshev's inequality we obtain, 
for $\IP$-almost all~$\omega$ and for all~$\vr\geq \vr_0$,
\begin{align}
 \Po^x[T^\vr>k] &\leq \Po^x[N^\vr_k(\Z_-)\geq k]\nonumber\\
            &\leq \Po^x[N^\vr_\infty(\Z_-)\geq k]\nonumber\\
          &\leq \frac{\Eo^{x}N^\vr_\infty(\Z_-)}{k}\nonumber\\
     &\leq \frac{C_1|x|}{k}.
\label{Cheb_tau}
\end{align}

Let us fix~$\delta_1>0$ such that $1+\delta_1<\alpha$ with $\alpha$ from Condition~C, 
and fix some $\beta\in\big(1,\frac{\alpha}{1+\delta_1}\big)$. 
Observe that~\eqref{Cheb_tau} implies that for any 
$x\in [-n^\beta,0]$
\begin{equation}
\label{oc_P_tau}
 \Po^x[T^\vr > n^{\beta(1+\delta_1)}] \leq C_1 n^{-\beta\delta_1}
\end{equation}
for $\IP$-almost all~$\omega$ and for all~$\vr\geq \vr_0$.

Fix a real number~$s\geq 1$ and denote
\[
 \sigma_s = \min\{k\geq 0 : S^\vr_k\in [-s,0]\}.
\]
Let~$G_s$ be  the event defined as
\[
 G_s=\{|S^\vr_{k}-S^\vr_{k-1}|\leq s \text{ for all }
 k\leq s^{\beta(1+\delta_1)}\}.
\]
By Condition~C, it is straightforward to obtain that, for some 
$C_2>0$
\begin{equation}
\label{est_G}
 \Po[G_s] \geq 1-C_2s^{-(\alpha-(1+\delta_1)\beta)}
\end{equation}
for all $\vr\geq \vr_0$
(observe that, by definition, $\alpha>(1+\delta_1)\beta$).
Also, note that, by the Markov property, when $x\leq y$ and $-s\leq y \leq 0$,
\begin{equation}
\label{eq_Mp}
 r_y^\vr = \Po^x[S^\vr_{T^\vr}=0\mid \sigma_s<T^\vr,
S^\vr_{\sigma_s}=y].
\end{equation}
On the event $G_s\cap\{T^\vr\leq s^{\beta(1+\delta_1)}\}$ we have 
$\sigma_s<T^\vr$ a.s. for $x<0$, and using 
also~\eqref{oc_P_tau}, \eqref{est_G},
\eqref{eq_Mp}, we obtain for any $x\in [-s^\beta,0)$
\begin{align}
 r_x^\vr &= \Po^x[S^\vr_{T^\vr}=0] \nonumber\\
   &\geq  \Po^x[S^\vr_{T^\vr}=0\mid
\sigma_s<T^\vr]\Po^x[\sigma_s<T^\vr]
\nonumber\\
  &\geq  \Big(\min_{y\in[-s,0]}r_y^\vr\Big) 
          \Po^x[G_s, T^\vr\leq s^{\beta(1+\delta_1)}] \nonumber\\
 &\geq  \Big(\min_{y\in[-s,0]}r_y^\vr\Big)
     (1-C_1s^{-\beta\delta_1}-C_2s^{-(\alpha-(1+\delta_1)\beta)}).
\label{conta_rx}
\end{align}

For any real number $k \geq 1$, define
\[
 u_k = \essinf_{\IP} \min_{y\in[-k,0]}r_y^\vr,
\]
which depends also on~$\omega$,
and let $\phi:=\min\{\beta\delta_1,\alpha-(1+\delta_1)\beta\}$.
Now, \eqref{conta_rx} implies that for some~$C_3>0$
\begin{equation}
\label{transfer}
 u_{s^\beta} \geq (1-C_3s^{-\phi})u_s.
\end{equation}
%%%%%% FC: need $n-0$ to start recursion
%% By the ellipticity Condition~E, we have~$u_2>0$. 
%%So, iterating~\eqref{transfer} we obtain that $u_m\geq 2\eps_1>0$
%% for all~$m\geq 2$, where
%%\[
%% \eps_1= \frac{1}{2} u_2(1-C_3 2^{-\phi})(1-C_3 2^{-\beta\phi})
%%    (1-C_3 2^{-\beta^2\phi}) (1-C_3 2^{-\beta^3\phi})\ldots
%%\]
%%(use the fact that $\sum_j 2^{-\beta^j\phi}<\infty$
%%to obtain that
%%the infinite product is strictly positive). 
By the ellipticity Condition~E, we have~$u_{2}>0$. 
Iterating~\eqref{transfer}, we obtain that $u_m\geq 2\eps_1>0$
 for all~$m\geq 2$, where
\[
\eps_1=\frac{1}{2} u_{2}(1-C_3 {2}^{-\phi})(1-C_3 {2}^{-\beta\phi})
    (1-C_3 {2}^{-\beta^2\phi}) (1-C_3 {2}^{-\beta^3\phi})\ldots
\]
is indeed positive since it holds that $\sum_j
{2}^{-\beta^j\phi}<\infty$.
This concludes the proof of Lemma~\ref{l_unif_eps}.
\qed

\medskip

Now, fix some integer $\vr\in [\vr_0,\infty)$, and consider 
a sequence of i.i.d.\ random variables
$\zeta_1,\zeta_2,\zeta_3,\ldots$ with
$P[\zeta_j=1]=1-P[\zeta_j=0]=\eps_1$ (the
parameter~$\eps_1$ is from Lemma~\ref{l_unif_eps}, and~$P$ stands 
for the law of this sequence; in the sequel we shall use also~$E$
for the expectation corresponding to~$P$). 
Then, our strategy can be
described in words in the following way. For all~$j\geq 1$,
Lemma~\ref{l_unif_eps} implies that $r_x^\vr(j\vr)\geq 2\eps_1$
for all $x\in[(j-1)\vr,j\vr-1]$. We couple the sequence
$\zeta=(\zeta_1,\zeta_2,\zeta_3,\ldots)$ with 
the random walk~$S^\vr$ in such a way that $\zeta_j=1$ implies
that $S^\vr_{T^\vr_{j\vr}}=j\vr$.
Denote 
\begin{equation}
\label{df_ell1}
 \ell_1=\min\{j : \zeta_j=1\}.
\end{equation}
Then, since~$\zeta$ (and therefore~$\ell_1$) is independent 
of~$\omega$, $\theta_{\ell_1\vr}\omega$ has the same law~$\IP$
as~$\omega$. This allows us to break the trajectory of 
the random walk into stationary ergodic (after suitable shift)
sequence of pieces, and then apply the 
ergodic theorem to obtain the
law of large numbers. 
The stationary measure of the environment seen
from the particle (for the truncated random walk) can also be
obtained from this construction by averaging along the cycle.
Then, we pass to the limit as $\rho\to\infty$.
\medskip

So, let us now construct the quenched law~$\Pozx$, i.e., 
the law of the random walk~$S^\vr$ when 
\emph{both} the environment~$\omega$ and
the sequence~$\zeta$ are fixed. This is done inductively: 
first, the law of $(S^\vr_k, k\leq T^\vr_\vr)$ is defined by 
\begin{align*}
 &\1{\zeta_1=1}\Po[\;\cdot\mid S^\vr_{T^\vr_\vr}=\vr]
 + \1{\zeta_1=0} \Big(\frac{r_0^\vr(\vr)-\eps_1}{1-\eps_1}
\Po[\;\cdot\mid S^\vr_{T^\vr_\vr}=\vr]\\
 &\hspace{7cm} {} + \frac{1-r_0^\vr(\vr)}{1-\eps_1}
\Po[\;\cdot\mid S^\vr_{T^\vr_\vr}>\vr] \Big).
\end{align*}
Then, given $S^\vr_{T^\vr_{j\vr}}=y\in[j\vr, (j+1)\vr-1]$, 
the law of $(S^\vr_k, T^\vr_{j\vr}+1\leq k\leq T^\vr_{(j+1)\vr})$ 
is
\begin{align*}
 &\1{\zeta_{j+1}=1}\Po^y[\;\cdot\mid
S^\vr_{T^\vr_{(j+1)\vr}}=(j+1)\vr]
\\
 & \hspace{1cm} + \1{\zeta_{j+1}=0} 
\Big(\frac{r_0^\vr((j+1)\vr)-\eps_1}{1-\eps_1}
          \Po^y[\;\cdot\mid S^\vr_{T^\vr_{(j+1)\vr}}=(j+1)\vr]\\
 &\hspace{5cm} {} + \frac{1-r_0^\vr((j+1)\vr)}{1-\eps_1}
         \Po[\;\cdot\mid S^\vr_{T^\vr_{(j+1)\vr}}>(j+1)\vr] \Big).
\end{align*}

Let $\IP':=\IP\otimes P$ (where~$P$ is the law of~$\zeta$), 
and~$\IE'$ the expectation corresponding to~$\IP'$.
With $\ell_0:=0$, let us define consistently with~\eqref{df_ell1}
\begin{equation}
\label{df_ellk}
 \ell_{k+1}=\min\{j>\ell_k : \zeta_j=1\}, \quad k\geq 0.
\end{equation}

Note that, by construction, 
\[
S^\vr_{T^\vr_{\ell_k\vr}}=\ell_k\vr
\]
 for all $k\geq 1$. 
We now define a regeneration structure, which is fundamental in our construction.
Following Chapter~8 of~\cite{Th}, we recall  the definition of
cycle-stationarity of a stochastic process together with a sequence
of points. Consider, on some probability space, (i) a sequence $Z=(Z_n)_n$ of random variables with
values in some measurable space  $(F, \mathcal F)$, (ii) a sequence of random times $\Sigma$ (called "time points"), $0<\Sigma_1<\Sigma_2<\ldots$. Define the $k$th cycle ${\mathcal C}_k=
(Z_n; \Sigma_k \leq n \leq \Sigma_{k+1}-1) \in \cup_{m \geq 1} F^m$. The sequence $Z$ is
cycle-stationary with points $\Sigma$ if $({\mathcal C}_k; k \geq 1)$ has the same law as
$({\mathcal C}_{k+1}; k \geq 1)$. It is cycle-stationary  and ergodic if $({\mathcal C}_k; k \geq 1)$ is stationary and ergodic.

\begin{lmm}
\label{l_stationary} Let $\rho < \8$.
The pair $(\theta_{S^\vr_\cdot}\omega, T^\vr_{\ell_\cdot\vr})$ 
is cycle-stationary and ergodic.  In particular,
$\theta_{\ell_k\vr}\omega$ has
the same law as~$\omega$ for all $k=1,2,3,\ldots$.
\end{lmm}
In short, the $k$th cycle ${\mathcal C}_k$ is the sequence of 
environments seen from the truncated walk $S^\vr_\cdot$ from 
time $T^\vr_{\ell_{k-1}\vr}$ to time $T^\vr_{\ell_{k}\vr}-1$ 
($k =1,2,\ldots$).
The first statement in the lemma is that the sequence 
$({\mathcal C}_k; k \geq 1)$ is
stationary under the measure $P \otimes \IP \otimes \Pozx$.
\smallskip

\noindent
\textit{Proof.} Let us denote by ${\bf {\mathcal C}}$ 
the above sequence,
and by $\vartheta$ the shift 
$(\vartheta {\bf {\mathcal C}})_k=
{\mathcal C}_{k+1}$. 
With $f \geq 0$ a measurable function on the appropriate space, 
we write
\begin{align*}
  E\IE \Eozx f(\vartheta {\bf {\mathcal C}}) 
&=
\sum_{m\geq 1}   E\IE \1{\ell_1=m}
\Eozx f(\vartheta {\bf {\mathcal C}})\\ 
&=
\sum_{m\geq 1}   E\IE \1{\ell_1=m}
\Eozx^{m \vr} f({\bf {\mathcal C}}) \qquad \qquad 
\text{(Markov\ property)}  \\ 
&=
\sum_{m\geq 1}  \eps_1(1-\eps_1)^{m-1}  \IE 
\Eozx^{m \vr} f({\bf {\mathcal C}}) \qquad \qquad
 \text{(independence)}  \\ 
&=
\sum_{m\geq 1}  \eps_1(1-\eps_1)^{m-1}  \IE 
\Eozx f({\bf {\mathcal C}}) \qquad \qquad 
\text{($\IP$-stationarity)}\\
&=  E\IE \Eozx f( {\bf {\mathcal C}}), 
\end{align*}
which shows the cycle-stationarity. The ergodicity then
follows from the ergodicity of~$\IP$ and independence of~$\omega$
and~$\zeta$, see Section~7 of Chapter~8 of~\cite{Th}. 
\qed

Now, we are able to prove the existence of the speed~$v_\vr$ 
for the truncated random walk.
First, we prove the following
\begin{lmm}
\label{l_expect_tau}
There exist $C_4,C_5>0$ such that for $\IP$-almost all~$\omega$ and
for all~$\vr \in [\vr_0, \8)$ we have
\begin{equation}
\label{eq_expect_tau}
C_4\vr \leq E\Eozx T^\vr_{\ell_1\vr}\leq C_5\vr.
\end{equation}
\end{lmm}

\noindent
\textit{Proof.}
We begin by proving the second inequality in~\eqref{eq_expect_tau}.
Write
%% FC:
%%\[
%% \Po^x[T^\vr>k] = \Po^x[N^\vr_k(\Z_-)\geq k] 
%%               \leq \Po^x[N^\vr_\infty(\Z_-)\geq k],
%%\]
%%so, using~\eqref{expected_hit_negative} we obtain
%%\begin{align}
%%\Eo^x T^\vr & \leq \sum_{k\geq 0}\Po^x[N^\vr_\infty(\Z_-)\geq k]
%%       \nonumber\\
%%      &= 1+\Eo^x N^\vr_\infty(\Z_-)\nonumber\\
%%   &\leq 1+C_1|x|.
%%\label{expect_tau}
%%\end{align}
%%%%%%%%%%%%%%%%%%
\[
 \Po^x[T^\vr>k] \leq \Po^x[N^\vr_k(\Z_-^*)\geq k+1] 
               \leq \Po^x[N^\vr_\infty(\Z_-)\geq k+1],
\]
so, using~\eqref{expected_hit_negative} we obtain
\begin{align}
\Eo^x T^\vr &= \sum_{k\geq 0}\Po^x[T^\vr>k] \nonumber\\
& \leq \sum_{k\geq 0}\Po^x[N^\vr_\infty(\Z_-)\geq k+1]
       \nonumber\\
      &= \Eo^x N^\vr_\infty(\Z_-)\nonumber\\
   &\leq C_1|x|.
\label{expect_tau}
\end{align}

Using the elementary inequality $\mathtt{E}(Y\mid A)\leq 
      \frac{\mathtt{E}Y}{\mathtt{P}[A]}$ for $Y \geq 0$
together with Lemma~\ref{l_unif_eps}, we obtain
that on~$\{\zeta_1=1\}$
\[
 \Eozx T^\vr_\vr 
    = \Eo \big(T^\vr_\vr\mid S^\vr_{T^\vr_\vr}=\vr\big) 
\leq \frac{\Eo T^\vr_\vr}{r^\vr_0(\vr)}
 \leq \frac{1}{\eps_1} \Eo T^\vr_\vr,
\]
and that on~$\{\zeta_1=0\}$ 
(observe that $\frac{r^\vr_0(\vr)-\eps_1}{r^\vr_0(\vr)}\leq 
\frac{1}{\eps_1}$ by Lemma~\ref{l_unif_eps})
\begin{align*}
 \Eozx T^\vr_\vr &= \frac{r^\vr_0(\vr)-\eps_1}{1-\eps_1} 
  \Eo \big(T^\vr_\vr\mid S^\vr_{T^\vr_\vr}=\vr\big) 
        + \frac{1-r^\vr_0(\vr)}{1-\eps_1} 
  \Eo \big(T^\vr_\vr\mid S^\vr_{T^\vr_\vr}>\vr\big)\\
 &= \frac{r^\vr_0(\vr)-\eps_1}{(1-\eps_1)r^\vr_0(\vr)} r^\vr_0(\vr)
  \Eo \big(T^\vr_\vr\mid S^\vr_{T^\vr_\vr}=\vr\big) \\
   &\qquad + \frac{1}{1-\eps_1} (1-r^\vr_0(\vr))
  \Eo \big(T^\vr_\vr\mid S^\vr_{T^\vr_\vr}>\vr\big)\\
 & \leq \frac{1}{\eps_1(1-\eps_1)}  
          \Eo \big(T^\vr_\vr
%\mid S^\vr_{T^\vr_\vr}>\vr
\big),
\end{align*}
so for any~$\zeta$ we obtain
\begin{equation}
\label{estim_tau_vr}
\Eozx T^\vr_\vr \leq \frac{1}{\eps_1(1-\eps_1)} \Eo T^\vr_\vr.
\end{equation}
In the same way, we show that for any~$\zeta$ and for all~$j\geq 1$
\begin{equation}
\label{estim_tau_others}
\Eozx^y T^\vr_{(j+1)\vr} \leq \frac{1}{\eps_1(1-\eps_1)} 
          \Eo^y T^\vr_{(j+1)\vr} 
\end{equation}
for all $y\in [j\vr, (j+1)\vr-1]$.

Writing
\[
 T^\vr_{k\vr} = T^\vr_\vr + (T^\vr_{2\vr}-T^\vr_\vr) 
         + \cdots + (T^\vr_{k\vr}-T^\vr_{(k-1)\vr}),
\]
and using~\eqref{expect_tau}, \eqref{estim_tau_vr}, 
\eqref{estim_tau_others},
we obtain on~$\{\ell_1=k\}$ that
\[
%% FC: 
%% \Eozx T^\vr_{k\vr} \leq \frac{1}{\eps_1(1-\eps_1)}k(1+C_1\vr).
 \Eozx T^\vr_{k\vr} \leq \frac{1}{\eps_1(1-\eps_1)}kC_1\vr.
\]
Since $P[\ell_1=k]=\eps_1(1-\eps_1)^{k-1}$, we see that, for
$\IP$-a.a.~$\omega$
%% FC: 
%%\begin{align*}
%% E\Eozx T^\vr_{\ell_1\vr} &\leq \sum_{k=1}^\infty
%%(1-\eps_1)^{k-2}k(1+C_1\vr)\\
%%&\leq C_5\vr
%%\end{align*}
%%for some $C_5>0$, and the proof of the second inequality
%%in~\eqref{eq_expect_tau} is finished.
\begin{align*}
 E\Eozx T^\vr_{\ell_1\vr} &\leq \sum_{k=1}^\infty
(1-\eps_1)^{k-2}kC_1\vr\\
= C_5\vr
\end{align*}
with $C_5=C_1 \eps_1^{-2}(1-\eps_1)^{-1}$,
and the proof of the second inequality
in~\eqref{eq_expect_tau} is finished.

Let us prove the first inequality in~\eqref{eq_expect_tau}. 
Consider a sequence of i.i.d.\ positive integer-valued 
random variables $Y_1,Y_2,Y_3,\ldots$ with the law
\[
 \mathtt{P}[Y_1 \geq s] = (\gamma_1s^{-\alpha}) \wedge 1,
\]
$s\geq 1$, with $\alpha, \gamma_1$ from Condition~C. Then, on 
$\{S^\vr_0=0\}$, 
it holds that for any~$\vr$ and for $\IP$-almost all~$\omega$,
$S^\vr_n$ is dominated by $Y_1+\cdots+Y_n$. From this, we easily
obtain
\begin{equation}
\label{lower_tau_omega}
 \Eo T^\vr_\vr \geq C_7\vr.
\end{equation}

Next, note that $r_0^\vr(\vr)-\eps_1\geq \eps_1$, 
so on the event~$\{\zeta_1=0\}$ we have
\begin{align}
 \Eozx T^\vr_\vr &= 
 \frac{r^\vr_0(\vr)-\eps_1}{(1-\eps_1)r^\vr_0(\vr)} 
r^\vr_0(\vr) \Eo \big(T^\vr_\vr\mid S^\vr_{T^\vr_\vr}=\vr\big) 
\nonumber\\
  &\qquad {}+ \frac{1}{1-\eps_1} (1-r^\vr_0(\vr))
  \Eo \big(T^\vr_\vr\mid S^\vr_{T^\vr_\vr}>\vr\big)\nonumber\\
 &\geq \frac{\eps_1}{1-\eps_1} \Eo T^\vr_\vr.
\label{*page9}
\end{align}
So, using~\eqref{lower_tau_omega} we obtain that on~$\{\zeta_1=0\}$
\[
 \Eozx T^\vr_{\ell_1\vr} \geq \Eozx T^\vr_\vr 
\geq \frac{\eps_1}{1-\eps_1} \Eo T^\vr_\vr \geq
\frac{C_7\eps_1}{1-\eps_1}\vr,
\]
and, finally,
\[
 \IE'\Eozx T^\vr_{\ell_1\vr} \geq \frac{C_7\eps_1}{1-\eps_1}\vr 
                P[\zeta_1=0] = C_7\eps_1\vr.
\]
Proof of Lemma~\ref{l_expect_tau} is finished.
\qed

Now, we show that
\begin{equation}
\label{eq_speed}
 v_\vr = \frac{\vr E\ell_1}{\IE'\Eozx T^\vr_{\ell_1\vr}}
        = \frac{\vr}{\eps_1\IE'\Eozx T^\vr_{\ell_1\vr}},
\end{equation}
which  implies that~$v_\vr>0$ by Lemma~\ref{l_expect_tau}. 
Indeed, suppose that~$n$ is such that
$T^\vr_{\ell_k\vr}\leq n < T^\vr_{\ell_{k+1}\vr}$. Then, we have
$\ell_{k+1}\vr - (T^\vr_{\ell_{k+1}\vr}-T^\vr_{\ell_k\vr})\vr <
S^\vr_n <
\ell_{k+1}\vr$, so
\begin{equation}
\label{eval_speed}
\frac{\ell_{k+1}\vr - (T^\vr_{\ell_{k+1}\vr}-T^\vr_{\ell_k\vr})\vr}
            {T^\vr_{\ell_{k+1}\vr}}
  \leq \frac{S^\vr_n}{n} \leq
\frac{\ell_{k+1}\vr}{T^\vr_{\ell_k\vr}}.
\end{equation}
Now, we divide the numerator and the denominator by~$k$ 
in~\eqref{eval_speed}, we use 
Lemmas~\ref{l_stationary},~\ref{l_expect_tau} and the
ergodic theorem to get
\[
\lim_{n \to \8}  \frac{S^\vr_n}{n} =v_\vr,
\]
with $v_\vr$ given by the second member of~\eqref{eq_speed}.
Since $(\ell_{k+1}-\ell_k)$ 
has a geometric distribution with parameter~$\eps_1$
we get the last expression in~\eqref{eq_speed}. 
This ends the proof of Theorem~\ref{t_rw_lln}
for finite $\vr$.

With the help of Lemma~\ref{l_stationary}, we derive that there 
exists an invariant measure for the environment seen from the
particle in the truncated case. By formula~(4.14${}^\circ$) of
Chapter~8 of~\cite{Th}, for $\vr<\infty$ we can characterize this
measure~$\IQ^\vr$ by its expectation  $\IE_{\vr}$,
\begin{equation}
\label{eq_Q}
 \IE_{\vr} f(\omega) = \frac{1}{\IE'\Eozx T^\vr_{\ell_1\vr}}
            \IE' \Eozx \sum_{k=1}^{T^\vr_{\ell_1\vr}}
                 f(\theta_{S^\vr_k}\omega).
\end{equation}
Next, we need to pass to the limit as $\vr\to\infty$. 
This requires a fine analysis of the Radon-Nikodym derivative
$\frac{d\IQ^{\vr}}{d\IP}$, and this is what we are going to do now.

\begin{prop}
\label{p_RN_bounded}
Let $\vr$ be finite. Then, 
\begin{equation}
\label{Rad_Nik}
 \frac{d\IQ^{\vr}}{d\IP}(\omega) 
= \frac{1}{\IE'\Eozx T^\vr_{\ell_1\vr}}
   \sum_{x\in\Z} E \mathtt{E}_{\theta_{-x}\omega^\vr,\zeta} 
N^\vr_{T^\vr_{\ell_1\vr}}(x).
\end{equation}
Moreover, there exist~${\hat\gamma}_1,{\hat\gamma}_2\in (0,\infty)$
(not depending on~$\vr$) such that for $\IP$-almost all~$\omega$
\begin{equation}
\label{eq_RN_bounded}
 {\hat\gamma}_1 \leq \frac{d\IQ^{\vr}}{d\IP}(\omega) \leq
{\hat\gamma}_2.
\end{equation}
\end{prop}

\noindent
\textit{Proof.} Using translation invariance of $\IP$, 
we write expression~\eqref{eq_Q} as

\begin{align*}
 \IE_{\vr} f(\omega) &= \frac{1}{\IE'\Eozx T^\vr_{\ell_1\vr}}
                             \int d\IP' \,       
\Eozx\sum_{k=1}^{T^\vr_{\ell_1\vr}}f(\theta_{S^\vr_k}\omega)\\
 &= \frac{1}{\IE'\Eozx T^\vr_{\ell_1\vr}} 
         \int d\IP'\sum_{k=1}^\infty 
              \Eozx(f(\theta_{S^\vr_k}\omega);k\leq
T^\vr_{\ell_1\vr})\\
 &= \frac{1}{\IE'\Eozx T^\vr_{\ell_1\vr}} 
  \sum_{k=1}^\infty \sum_{x\in\Z} \int d\IP'\,
                       f(\theta_x\omega) \Pozx[S^\vr_k=x,k\leq
T^\vr_{\ell_1\vr}]\\
  &= \frac{1}{\IE'\Eozx T^\vr_{\ell_1\vr}} 
  \sum_{k=1}^\infty \sum_{x\in\Z} \int d\IP' \, f(\omega) 
    \mathtt{P}_{\!\theta_{-x}\omega^\vr,\zeta}
                    [S^\vr_k=x,k\leq T^\vr_{\ell_1\vr}]\\
 &= \int d\IP' \, f(\omega) \frac{1}{\IE'\Eozx T^\vr_{\ell_1\vr}} 
    \sum_{x\in\Z} \mathtt{E}_{\theta_{-x}\omega^\vr,\zeta}
N^\vr_{T^\vr_{\ell_1\vr}}(x)\\
 &= \int d\IP \, f(\omega) \frac{1}{\IE'\Eozx T^\vr_{\ell_1\vr}}
\sum_{x\in\Z} 
          E\mathtt{E}_{\theta_{-x}\omega^\vr,\zeta}
                 N^\vr_{T^\vr_{\ell_1\vr}}(x),
\end{align*}
which proves~\eqref{Rad_Nik}.

Let us prove~\eqref{eq_RN_bounded}. Write
\begin{align}
 \frac{d\IQ^{\vr}}{d\IP}(\omega) &= \frac{1}{\IE'\Eozx 
T^\vr_{\ell_1\vr}}
E\Big(\sum_{x<0}\mathtt{E}_{\theta_{-x}\omega^\vr,\zeta}
N^\vr_{T^\vr_{\ell_1\vr}}(x)
    +\sum_{x\geq 0}\mathtt{E}_{\theta_{-x}\omega^\vr,\zeta}
N^\vr_{T^\vr_{\ell_1\vr}}(x)\Big)\nonumber\\
  &\leq \frac{1}{\IE'\Eozx T^\vr_{\ell_1\vr}} E \sum_{x<0}
 \mathtt{E}_{\theta_{-x}\omega^\vr,\zeta} 
         N^\vr_\infty(x)\nonumber\\
&\qquad +\frac{1}{\IE'\Eozx T^\vr_{\ell_1\vr}}
E \sum_{0<x\leq (\ell_1+1)\vr}
\mathtt{E}_{\theta_{-x}\omega^\vr,\zeta} N^\vr_\infty(x)\nonumber\\
 & =: \mathbf{A}_1 + \mathbf{A}_2.
\label{bound_on_RN}
\end{align}
Next, we need to obtain upper bounds on the terms
 $\mathbf{A}_1, \mathbf{A}_2$; for that, let us write first
\begin{equation}
\label{sum_N}
 N^\vr_\infty(x) =
N^\vr_{[0,T^\vr_\vr)}(x)+N^\vr_{[T^\vr_\vr,T^\vr_{2\vr})}(x)
                   +N^\vr_{[T^\vr_{2\vr},T^\vr_{3\vr})}(x)+\cdots.
\end{equation}
Analogously to~\eqref{estim_tau_vr} and~\eqref{estim_tau_others},
 for any~$\zeta$ we obtain
\begin{equation}
\label{estim_N_1}
 \Eozx N^\vr_{[0,T^\vr_\vr)}(x) \leq \frac{1}{\eps_1(1-\eps_1)}  
          \Eo N^\vr_{[0,T^\vr_\vr)}(x),
\end{equation}
and for all~$j\geq 1$ and all $y\in[j\vr, (j+1)\vr-1]$
\begin{equation}
\label{estim_N_others}
 \Eozx^y N^\vr_{[T^\vr_{j\vr},T^\vr_{(j+1)\vr})}(x) \leq \frac{1}
 {\eps_1(1-\eps_1)} 
   \Eo^y N^\vr_{[T^\vr_{j\vr},T^\vr_{(j+1)\vr})}(x).
\end{equation}

Consider the term~$\mathbf{A}_1$ of~\eqref{bound_on_RN}. 
Applying~\eqref{estim_N_1} and~\eqref{estim_N_others}
to~\eqref{sum_N} and using Condition~D, we obtain for $\IP$-almost 
every~$\omega$ (and so the following holds also with
$\theta_{-x}\omega^\vr$ on
the place of $\omega^\vr$) and all~$x<0$ that
\begin{align}
 \Eozx N^\vr_\infty(x) &\leq \frac{1}{\eps_1(1-\eps_1)}
\left(g_1(|x|) + g_1(\vr+|x|) + g_1(2\vr +
|x|) 
+ \cdots \right) \nonumber\\
 &\leq \frac{1}{\eps_1(1-\eps_1)}\sum_{j=0}^\infty g_1(|x|+j).
\label{est_EN(-x)}
\end{align}
So, we write
\begin{align}
 \mathbf{A}_1 &\leq \frac{1}{\eps_1(1-\eps_1)\IE'\Eozx 
 T^\vr_{\ell_1\vr}} \sum_{x<0} \sum_{j=0}^\infty
g_1(|x|+j)\nonumber\\
 &\leq \frac{1}{\eps_1(1-\eps_1)\IE'\Eozx T^\vr_{\ell_1\vr}} 
   \sum_{k=1}^\infty k g_1(k)\nonumber\\
 &\leq C_5
\label{bound_T1}
\end{align}
for some $C_5>0$.

Let us deal now with the term~$\mathbf{A}_2$. 
Suppose that~$x\geq 0$ is such that $x\in[k\vr,(k+1)\vr)$. Then, we
have $N^\vr_{[T^\vr_{j\vr},T^\vr_{(j+1)\vr})}(x)=0$ for all $j<k$.
So, by~\eqref{estim_N_others}, we obtain for $\IP$-almost
every~$\omega$ (again, this means that it holds also with
$\theta_{-x}\omega^\vr$ on the place of $\omega^\vr$) that
\begin{equation}
\label{N_near}
 \Eozx N^\vr_{[T^\vr_{j\vr},T^\vr_{(j+1)\vr})}(x) 
   \leq \frac{g_1(0)}{\eps_1(1-\eps_1)}
\end{equation}
for $j\in\{k,k+1\}$ 
%%%%(again, the function~$g_1$ is from Condition~D)
, and
\begin{equation}
\label{N_far}
 \Eozx N^\vr_{[T^\vr_{j\vr},T^\vr_{(j+1)\vr})}(x) 
     \leq \frac{1}{\eps_1(1-\eps_1)} g_1((j-k-1)\vr)
\end{equation}
for $j>k+1$. Using~\eqref{sum_N}, we obtain 
\[
 \Eozx N^\vr_\infty (x) \leq \frac{2g_1(0)}{\eps_1(1-\eps_1)} 
          + \frac{1}{\eps_1(1-\eps_1)}\sum_{i\geq 1}g_1(i) 
      \leq C_6
\]
for some $C_6>0$. Then, using also Lemma~\ref{l_expect_tau}, 
we have
\begin{align}
 \mathbf{A}_2 &\leq \frac{1}{\IE'\Eozx T^\vr_{\ell_1\vr}}C_6\vr 
E (\ell_1+1)\nonumber\\
 &\leq \frac{C_6}{C_4} (\eps_1^{-1}+1).
\label{bound_T2}
\end{align}
Using~\eqref{bound_T1} and~\eqref{bound_T2}, we obtain the second
inequality in~\eqref{eq_RN_bounded}.

As for the first inequality in~\eqref{eq_RN_bounded}, let us note
that, analogously to~\eqref{*page9}, on the event $\{\zeta_1=0\}$ 
we have for any nonnegative random variable~$Y$ which is measurable
with respect to the sigma-algebra generated by
$(S^\vr_1,\ldots,S^\vr_{T^\vr_\vr})$
\begin{equation}
\label{expect_Y}
 \Eozx Y \geq \frac{\eps_1}{1-\eps_1} \Eo Y.
\end{equation}
By Lemma~\ref{l_unif_eps}, we obtain
\begin{align*}
 E\sum_{x\in \Z}
\mathtt{E}_{\theta_{-x}\omega^\vr,\zeta}
     N^\vr_{T^\vr_{\ell_1\vr}}(x)
  &\geq (1-\eps_1) E\sum_{x\in [1,\vr]} 
\mathtt{E}_{\theta_{-x}\omega^\vr,\zeta}N^\vr_{T^\vr_{\ell_1\vr}}
(x)\\
 &\geq (1-\eps_1) \times\frac{\eps_1}{1-\eps_1} 
E\sum_{x\in [1,\vr]} 
\mathtt{E}_{\theta_{-x}\omega^\vr}N^\vr_{T^\vr_{\ell_1\vr}}(x)\\
&\geq 2\eps_1^2 \vr,
\end{align*}
and this concludes the proof of Proposition~\ref{p_RN_bounded}.
\qed

\medskip

Now, we pass to the limit as~$\vr\to\infty$, proving 
the existence of~$v^\infty,\IQ:=\IQ^\infty$, and the fact that
$v^\vr\to v^\infty$.

From Proposition~\ref{p_RN_bounded} we obtain that the family
of measures $(\IQ^\vr, \vr\in[\vr_0,+\infty))$ is tight. 
By Prohorov's theorem, 
there exists a sequence $\vr_k\to\infty$ and a
probability measure~$\IQ^\infty$ such that $\IQ^{\vr_k}\to
\IQ^\infty$ weakly as $k\to\infty$. Let us prove that~$\IQ^\infty$ 
is in fact a stationary measure for the environment of the random
 walk without
truncation $S^\infty$. For this, take a bounded continuous
function~$f$ and, recalling that $\IE_\vr$ is the expectation
corresponding to~$\IQ^\vr$, write by stationarity
\begin{align}
 |\IE_\infty f(\omega) - \IE_\infty
f(\theta_{S_1^\infty}\omega)|\leq & |\IE_\infty f(\omega) -
\IE_{\vr_k} f(\omega)| \nonumber\\
 & {} + |\IE_{\vr_k} f(\theta_{S_1^{\vr_k}}\omega) -
\IE_\infty f(\theta_{S_1^{\vr_k}}\omega)| \nonumber\\
 & {} + \IE_\infty
|f(\theta_{S_1^{\vr_k}}\omega)-f(\theta_{S_1^\infty}\omega)|
\nonumber\\
=: & \mathbf{B}_1 + \mathbf{B}_2 + \mathbf{B}_3.
\label{conta_weak_limit}
\end{align}
First, we deal with terms~$\mathbf{B}_2$ and~$\mathbf{B}_1$.
By Condition~C, for any~$\eps>0$ there exists~$h_0$ such that for
any~$\vr$ we have
\[
 \Po[|S_1^\vr| > h_0] < \eps \qquad \text{$\IP$-a.s.}
\]
Then, supposing without restriction of generality that $|f|\leq 1$,
we write
\begin{align*}
 \mathbf{B}_2 &\leq \Big|\IE_{\vr_k}
\sum_{|m|\leq h_0}\omega_{0m}^{\vr_k}f(\theta_m\omega) - 
\IE_\infty
\sum_{|m|\leq h_0}\omega_{0m}f(\theta_m\omega)\Big|+2\eps\\
&\leq \Big|\IE_{\vr_k}
\sum_{|m|\leq h_0}\omega_{0m}f(\theta_m\omega) - 
\IE_\infty
\sum_{|m|\leq h_0}\omega_{0m}f(\theta_m\omega)\Big|
 + \IE_{\vr_k}\sum_{|m|>\vr_k}\omega_{0m} +2\eps.
\end{align*}
Since $\IQ^{\vr_k}$ converges to $\IQ^\infty$ as $k\to\infty$
and using also Condition~C, 
we can choose~$k$ large enough in such a way that
$\mathbf{B}_1+\mathbf{B}_2 \leq 5\eps$. As for the
term~$\mathbf{B}_3$, again using Condition~C we note 
that, for the natural coupling of $S_1^{\vr_k}$ and $S_1^\infty$
(described in Section~\ref{s_results_rwre}), we
have 
\[
 \Po[S_1^{\vr_k}\neq S_1^\infty] < \eps  \qquad \text{$\IP$-a.s.}
\]
for large enough~$k$. Then, for such~$k$'s, we have 
$\mathbf{B}_3<\eps$ (recall that we assumed that $|f|\leq 1$).
Then, since~$\eps$ is arbitrary, \eqref{conta_weak_limit} implies
that~$\IQ^\infty$ is stationary for~$S^\infty$. 

Now, let us prove that~$\IQ^\vr$ is ergodic for
all~$\vr\in[\vr_0,+\infty]$. First, we note
that~\eqref{eq_RN_bounded} holds for~$\vr=\infty$ as well. Then, we
argue by contradiction: suppose that~$\IQ^\vr$ is not ergodic and
let~$A\subset\Omega$ be a nontrivial invariant event.
Then, \eqref{eq_RN_bounded} implies that $0<\IP[A]<1$. 
Since $\IP$ is stationary
ergodic under space shift, the random set
\[
{\mathfrak G}(\omega)=\{k>0: \theta_k\omega\in A\}
\]
is such that 
\[
|{\mathfrak G}|=|{\mathfrak G}^c|=\infty \qquad \text{$\IP$-a.s.},
\]
which means by absolute continuity that
\[
|{\mathfrak G}|=|{\mathfrak G}^c|=\infty \qquad
\text{$\IQ^\vr$-a.s.}
\]
Since they are infinite, it follows from Lemma~\ref{l_unif_eps} 
that both sets ${\mathfrak G}$
and ${\mathfrak G}^c$ are visited by~$S^\vr$
infinitely many times almost surely. 
This contradicts the invariance of $A$ 
for the environment dynamics, which asserts that 
$\theta_{S_k}\omega \in A$ holds
$\IQ^\vr \otimes \Po$-a.s.\ on $\{\omega \in A\}$ and that 
$\theta_{S_k}\omega \in A^c$ holds
$\IQ^\vr \otimes \Po$-a.s.\ on $\{\omega \in A^c\}$.

Now, we claim that the ergodicity of $\IQ^\infty$ implies that
$\IQ^\vr \to \IQ^\infty$ weakly as $\vr\to\infty$: indeed, any 
possible limit~$\tilde\IQ$ has also to be
ergodic by the above argument, and thus singular with respect to~$\IQ^\infty$ or equal 
to it; in view of~\eqref{eq_RN_bounded},
the first case cannot hold, so the second one does hold.

Finally, let us deal with the proof 
of~\eqref{eq_rw_speed} and the fact
that $v_\vr\to v_\infty$ as $\vr\to\infty$. 
The existence of~$v_\vr$ (which has been previously established 
for finite~$\vr$) and
formula~\eqref{eq_rw_speed} follow from the ergodic theorem. To
prove the convergence, note that, by~\eqref{eq_rw_speed},
\[
 v_\vr = \IE_\vr \sum_{m\in\Z}m\omega_{0m},
\]
and, by Condition~C, for any~$\eps>0$
there exists~$h_1$ such that 
\begin{equation}
\label{hvosty}
\sum_{|m|>h_1}|m|\omega_{0m} < \eps \qquad \text{$\IP$-a.s.}
\end{equation}
Since~$\eps$ is arbitrary, 
the fact that $v_\vr\to v_\infty$ follows from~\eqref{hvosty} and
 from
\[
 \IE_\vr\sum_{|m|\leq h_1}m\omega_{0m} \to
\IE_\infty\sum_{|m|\leq h_1}m\omega_{0m} \qquad \text{as
$\vr\to\infty$}.
\]
This concludes the proof of Theorems~\ref{t_rw_lln}
and~\ref{t_rw_env}. \qed

\section{Proofs for the random billiard}
\label{s_proofs_billiard}
 Let us define $h_0:=(\eps^{-1}\gamma_d)^{1/(d-1)}$, with~$\eps$
  from Condition~P and $\gamma_d$ from (\ref{def_K}). From~\eqref{def_K} it directly follows that
\begin{equation}
\label{bound_K}
 K(x,y) \leq \frac{\gamma_d}{\|x-y\|^{d-1}}.
\end{equation}
So, if Condition~P holds, for all pairs of points $z,z'$ 
involved in Condition~P it holds that $\|z-z'\|^{d-1}\leq
\eps^{-1}\gamma_d$, so that
\[
 |(z-z')\cdot\e| \leq \|z-z'\| \leq h_0. 
%% (\eps^{-1}\gamma_d)^{1/(d-1)}.
\]
Thus, we obtain that Condition~P holds
for~${\hat K}$ on the place of~$K$ with~$\eps'=\eps
e^{-\lambda h_0}$ on the place of~$\eps$. 

% Define
% \[
%  {\hat F}^\omega(x) = \{y\in\partial\omega : -1\leq (y-x)\cdot\e 
% \leq 1\}.
% \]
Let us consider a sequence of i.i.d.\ random
variables~$\eta_1,\eta_2,\eta_3,\ldots$ with uniform 
distribution on
$\{1,2,\ldots,N\}$ (where~$N$ is from Condition~P), independent of
everything else. Also, 
define $J(0)=0$, $J(n)=\eta_1+\cdots+\eta_n$.
Then, analogously to the proof of Lemma~3.6 in~\cite{CPSV2},
using Condition~P one can obtain that, for
any~$x\in\partial\omega$ 
and~$B\subset \{y\in\partial\omega : -1\leq
(y-x)\cdot\e \leq 1\}$, we have
\begin{equation}
\label{oc_plotnost'}
 \Po^x[\xi_{\eta_1}\in B] \geq \frac{1}{N}\delta^{N-1}(\eps
e^{-\lambda h_0})^N  \nu^\omega(B).
\end{equation}
Here, and below, we still use, for simplicity, the notations 
$\Po^x, \Eo^x$ for the enlarged model, meaning that
$\omega$ is fixed but $\xi_\cdot$ and $\eta_\cdot$ 
are integrated out.

% Consistently with the notations of~\cite{CPSV2} (see the proof
% of Proposition~4.1 there), we 
Let us define for an arbitrary $a,b\in\R$, $a<b$,
\[
{\tilde F}^\omega(a,b) = \{x\in\partial\omega : x\cdot
\e\in[a,b]\}.
\]
Then, consider the moment when the particle 
steps out from~${\tilde F}^\omega(a,b)$:
\[
 \tau(a,b) = \min\{n\geq 0: |\xi_n\cdot \e| \notin [a,b]\}
  = \min\{n\geq 0: \xi_n \notin {\tilde F}^\omega(a,b)\}.
\]
Next, we prove the following fact:
\begin{lmm}
\label{l_go_out} 
There exist ${\tilde \gamma}_1,{\tilde \gamma}'_1>0$ such that
for $\IP$-almost all~$\omega$, 
for any $a\leq b-1$ we have
\begin{equation}
\label{tail_tau}
 \Po^x[\tau(a,b)>(b-a)^3t] \leq {\tilde \gamma}_1e^{-{\tilde
\gamma}'_1 t}
\end{equation}
for all~$x\in {\tilde F}^\omega(a,b)$ and all~$t\geq 1$.
\end{lmm}

\noindent
\textit{Proof.}
Observe that, from Condition~L we obtain that 
for some positive constants~${\tilde\gamma}_2,{\tilde\gamma}_3$
we have
\begin{equation}
\label{bounded_nuo}
  {\tilde\gamma}_2 \leq \nuo\big(x:x\cdot\e\in[s,s+1]\big)\leq
{\tilde\gamma}_3 \qquad \text{$\IP$-a.s.}
\end{equation}
for all~$s\in \R$ (without restriction of generality one may assume
that ${\tilde\gamma}_2\leq 1$, ${\tilde\gamma}_3\geq 1$).

Now, suppose without restriction of generality that $a=0$, and~$b$
is a positive integer. Denote 
\[
 {\tilde\tau} = \min\{n\geq 0: \xi_{J(n)} \notin 
{\tilde F}^\omega(a,b)\}.
\]
Clearly, by definition of the random variables $(\eta_i)$, we have
$\{\tau(a,b)>b^3t\}\subset \{{\tilde\tau}>N^{-1}b^3t\}$. 
Now, let us
consider a sub-Markov kernel $Q_\omega:=Q_\omega^1$, which acts on
functions $f:\partial\omega \to \R$ in the following way:
\begin{equation}
\label{def_Q}
(Q_\omega^n f)(x) = \Eo^x\big(f(\xi_{J(n)})\1{{\tilde\tau}>n}\big).
\end{equation}

Let
\[
 {\tilde K}(x,y) = \frac{1}{N}\sum_{j=1}^N {\hat K}^j(x,y)
\]
be the transition density of the process $(\xi_{J(n)}, n\geq 0)$.
Observe that this process is reversible with the reversible
measure~$\nuo_\lambda$, so that 
$\pi(x){\tilde K}(x,y)=\pi(y){\tilde K}(y,x)$ for all 
$x,y\in\partial\omega$. Let 
\[
 \EE(f,g) = \intl_{(\partial\omega)^2}\pi(x){\tilde K}(x,y)
  (f(x)-f(y))(g(x)-g(y)) \, d\nuo(x) \, d\nuo(y),
\]
and define 
\begin{equation}
\label{def_spectral_gap}
 {\tilde\Lambda} = \inf\Big\{\frac{\EE(f,f)}{2\int_{\partial\omega}
f^2(x)\, d\nuo_\lambda(x)} : f\not\equiv 0, 
f\big|_{({\tilde F}^\omega(0,b))^c}=0,f\in L^2(\nuo_\lambda)\Big\}.
\end{equation}
From the variational formula for the top of the 
spectrum of symmetric operators, 
 $\|Q_\omega\|_{L^2(\nuo_\lambda)}=1-{\tilde\Lambda}$,
so we look for  a lower bound for~${\tilde\Lambda}$. Denote
\[
 U_j = \{x\in\partial\omega : x\cdot\e\in (j,j+1]\}, \qquad 
  j\in\Z,
\]
so that ${\tilde F}^\omega(j,j+1)=\overline{U_j}$.
Observe that, by~\eqref{oc_plotnost'}, for any~$j\in\Z$ we have
\[
 {\tilde K}(x,y) \geq \frac{1}{N}\delta^{N-1}(\eps
e^{-\lambda h_0})^N
\]
for all $x\in U_j$, $y\in U_{j+1}$. Also, it is clear that
\[
\nuo(\{x\in\partial\omega: x\cdot \e = 0\})=0\quad \IP\text{-a.s.},
\]
and
\begin{equation}
\label{bound_pi_U}
 e^{\lambda n} \leq \pi(x) \leq e^{\lambda (n+1)}
 \qquad \text{for all }x\in U_n.
\end{equation}
So, using also~\eqref{bounded_nuo}, the
Cauchy-Schwarz inequality, and the fact that $f(y)=0$ for all~$y\in
U_b$, we can write
\begin{align*}
 \lefteqn{\intl_{\partial\omega} f^2(x)\, d\nuo_\lambda(x)}\\
 &=
\intl_{\partial\omega} \pi(x)f^2(x)\, d\nuo(x)\\
&= \sum_{i=0}^{b-1} \intl_{U_i} \pi(x)f^2(x)\, d\nuo(x)\\
 &= \sum_{i=0}^{b-1}\frac{1}{\nuo(U_{i+1})\ldots \nuo(U_b)}
 \intl_{U_i}\pi(x_i) \, d\nuo(x_i)\intl_{U_{i+1}}d\nuo(x_{i+1})\\
 & \qquad\qquad\qquad  \ldots
\intl_{U_b}d\nuo(x_b)\Big(\sum_{j=i}^{b-1}
  (f(x_j)-f(x_{j+1}))\Big)^2\\
&\leq b \sum_{i=0}^{b-1}\frac{1}{\nuo(U_{i+1})\ldots \nuo(U_b)}
 \intl_{U_i}\pi(x_i) \, d\nuo(x_i)\intl_{U_{i+1}}d\nuo(x_{i+1})\\
 & \qquad\qquad\qquad  \ldots
\intl_{U_b}d\nuo(x_b)\Big(\sum_{j=i}^{b-1}
  (f(x_j)-f(x_{j+1}))^2\Big)\\
&\leq b\sum_{i=0}^{b-1} e^{\lambda(i+1)}
\Big(\frac{1}{\nuo(U_{i+1})}
 \intl_{U_i}d\nuo(x_i)\intl_{U_{i+1}}d\nuo(x_{i+1})
  (f(x_i)-f(x_{i+1}))^2 \\
 & \qquad ~~~ +\sum_{j=i+1}^{b-1}
 \frac{\nuo(U_i)}{\nuo(U_j)\nuo(U_{j+1})}
  \intl_{U_j}d\nuo(x_j)\intl_{U_{j+1}}d\nuo(x_{j+1})
  (f(x_j)-f(x_{j+1}))^2\Big)\\
 & \leq \frac{b {\tilde\gamma}_3}{{\tilde\gamma}_2^2}
 \sum_{i=0}^{b-1} e^{\lambda(i+1)} \sum_{j=i}^{b-1}
  \intl_{U_j}d\nuo(x_j)\intl_{U_{j+1}}d\nuo(x_{j+1})
  (f(x_j)-f(x_{j+1}))^2
   \quad \text{(by~\eqref{bounded_nuo})} \\
 &\leq \frac{b {\tilde\gamma}_3}{{\tilde\gamma}_2^2} 
 \left(\frac{1}{N}\delta^{N-1}(\eps
e^{-\lambda h_0})^N\right)^{-1}\;\sum_{j=0}^{b-1}\Big(\sum_{i=0}^j
e^{\lambda(i+1)}\Big) e^{-\lambda j} \\
 & \qquad\qquad\times\intl_{U_j}d\nuo(x_j)
\intl_{U_{j+1}}d\nuo(x_{j+1}) \pi(x_j) {\tilde K}(x_j,x_{j+1})
  (f(x_j)-f(x_{j+1}))^2\\
 & \leq \frac{b {\tilde\gamma}_3 N }{{\tilde\gamma}_2^2  
 \delta^{N-1}(\eps
e^{-\lambda h_0})^N} \Big(\sum_{i=-\infty}^0
e^{\lambda(i+1)}\Big) \EE(f,f),
\end{align*}
and so, for some positive constant~$C_1$ it holds that 
\begin{equation}
\label{oc_spectral_gap}
 {\tilde\Lambda} \geq \frac{C_1}{b}.
\end{equation}

Then, since ${\tilde\Lambda} = 1-\|Q_\omega\|_{L^2(\nuo_\lambda)}$
from the spectral variational formula we have 
\begin{equation}
\label{L2_bound_Q}
\|Q_\omega^n\|_{L^2(\nuo_\lambda)}\leq \Big(1-\frac{C_1}{b}\Big)^n.
\end{equation}

Now, using the notation $\mathbf{1}_B$ 
for the indicator function of
$B\subset\partial\omega$, observe that 
$\Po^x[{\tilde\tau}>n]=(Q^n_\omega\mathbf{1}_{{\tilde
F}^\omega(0,b)})(x)$.
For $j\in[0,b-1]$ we can write
using~\eqref{bounded_nuo}, \eqref{L2_bound_Q}, and Cauchy-Schwarz
inequality
\begin{align*}
 \intl_{U_j} \Po^x[{\tilde\tau}>n] \, d\nuo_\lambda(x)
  &= (\mathbf{1}_{U_j},Q_\omega^n\mathbf{1}_{{\tilde
F}^\omega(0,b)})_{L^2(\nuo_\lambda)}\\
 &\leq \|\mathbf{1}_{U_j}\|_{L^2(\nuo_\lambda)} 
\|\mathbf{1}_{{\tilde F}^\omega(0,b)}\|_{L^2(\nuo_\lambda)}
 \|Q_\omega^n\|_{L^2(\nuo_\lambda)} \\
 & \leq C_2 e^{\lambda j/2}e^{\lambda b/2}
\Big(1-\frac{C_1}{b}\Big)^n
\end{align*}
for some~$C_2>0$. So, again using~\eqref{bound_pi_U}, 
we have for $j \leq b-1$,
\begin{equation}
\label{bound_int_Uj}
 \intl_{U_j} \Po^x[{\tilde\tau}>n] \, d\nuo(x)
 \leq C_2 e^{\lambda(b-j)/2}\Big(1-\frac{C_1}{b}\Big)^n.
\end{equation}
Now, \eqref{bound_int_Uj} implies that, if~$b$ is large enough, 
then
\[
 \intl_{U_j} \Po^x[{\tilde\tau}\leq b^3-1] \, d\nuo(x) \geq C_3>0.
\]
So, with $C_4:=C_3N^{-1}\delta^{N-1}(\eps
e^{-\lambda h_0})^N>0$, for any~$x\in U_j$ we can write 
(using also~\eqref{oc_plotnost'}) that 
\begin{align*}
 \Po^x[{\tilde\tau}\leq b^3] &=\intl_{\partial\omega}
  {\tilde K}(x,y) \Po^y[{\tilde\tau}\leq b^3-1]\, d\nuo(y)\\
& \geq \frac{1}{N}\delta^{N-1}(\eps
e^{-\lambda h_0})^N \intl_{U_j} \Po^y[{\tilde\tau}\leq b^3-1] \,
d\nuo(y) \\
& \geq C_4.
\end{align*}
This implies that $\Po^x[{\tilde\tau}>b^3t]\leq e^{-C_6t}$ 
for any $x\in{\tilde F}^\omega(0,b)$, and this (as discussed in the
beginning of the proof of this lemma) by its turn
implies~\eqref{tail_tau}, thus concluding the proof of
Lemma~\ref{l_go_out}.
\qed

\medskip

Consider~$B\subset\partial\omega$ with positive
 $(d-1)$-dimensional Hausdorff measure and such that
$\sup\{x\cdot \e : x\in B\}<+\infty$. By definition, the
\emph{stationary
distribution~$\pi^B$ conditioned on}~$B$ is the distribution with the density
$\frac{\pi(x)}{\pi(B)} {\bf 1}_{B}(x)$ 
(recall that $\pi(B):=\nuo_\lambda(B)$).
We use the notation $\Po^B$, $\Eo^B$ for the KRWD starting from the
stationary distribution conditioned on~$B$.

%%% FC: changed and added
% 
% Recall that $\eta_1,\eta_2,\eta_3,\ldots$ are i.i.d.\
% random variables with uniform distribution on $\{1,\ldots,N\}$, 
% and
%  $J(0)=0$, $J(n)=\eta_1+\cdots+\eta_n = J(n-1)+\eta_n$. 
% %% FC: suppressed def. here, and introduced later when needed
% %Let $\mathcal{F}_m$ be the sigma-algebra generated 
% by $\xi_0,\ldots,\xi_m$.
% We now focus on the process $(\xi_{J(n)}, n\geq 0)$. 
% Consider also a
% sequence of i.i.d.\ random variables $\zeta'_n\in\{0,1\}$ with 
% \[
%  P[\zeta'_n=1] = r_1
% \]
% where $r_1=\frac{1}{Ne^\lambda}{\tilde\gamma}_2 \delta^{N-1}(\eps
% e^{-\lambda h_0})^N$
% (recall~\eqref{oc_plotnost'} and~\eqref{bounded_nuo}). Then,
% one can couple the random sequence
% $(\xi_{J(n)},n\geq 1)$ with $\zeta'=(\zeta'_n,n\geq 1)$ in such a
% way that when the event $\{\zeta'_n=1\}$ occurs, $\xi_{J(n)}$ has
% the stationary distribution on $U_{[\xi_{J(n-1)}\cdot\e]}$.
% %%%%FC: tried to make it more clear 7 lines below.
% Check if it is OK
% %; 
% % moreover, conditionally on~$\xi_{J(n-1)}$, $\zeta'_n$
% % depends only on $(\eta_n,\xi_{J(n-1)+1},\ldots,\xi_{J(n)})$.
% We denote by~$\Pozz$ and~$\Eozz$ the probability and expectation
% with fixed~$\omega$ and~$\zeta'$, and let~$E^{\zeta'}$ be the
% expectation with respect to~$\zeta'$. 

Now, we construct the connection with RWRE. 
Recall that $\eta_1,\eta_2,\eta_3,\ldots$ are i.i.d.\
random variables with uniform distribution on $\{1,\ldots,N\}$, and
 $J(0)=0$, $J(n)=\eta_1+\cdots+\eta_n = J(n-1)+\eta_n$. 
We now focus on the process $(\xi_{J(n)}, n\geq 0)$. 
In view of~\eqref{oc_plotnost'} and~\eqref{bounded_nuo},
we couple this process with a Bernoulli process 
$\zeta'=(\zeta'_n, n\geq 1)$ (independent of $\omega$)
of parameter $r_1=(Ne^\lambda)^{-1}{\tilde\gamma}_2 
\delta^{N-1}(\eps e^{-\lambda h_0})^N$,
\[
 P[\zeta'_n=1] = 1-P[\zeta'_n=0] = r_1,
\]
in such a way that on the event $\{\zeta'_n=1\}$, $\xi_{J(n)}$ has
the stationary distribution 
on $U_{[\xi_{J(n-1)}\cdot\e]}$.
%%%%%%%%%%%%%  FC: added 30/08/10
(The choice of the stationary distribution is arbitrary, and the whole 
construction could be implemented with another distribution, 
absolutely continuous to the uniform on $U_{[\xi_{J(n-1)}\cdot\e]}$ with
density bounded from above and below.)
We denote by~$E^{\zeta'}$ be the
expectation with respect to~$\zeta'$, and~$\Pozz$,~$\Eozz$ 
the probability and expectation
with fixed~$\omega$ and~$\zeta'$, which is defined as follows: 
$\Pozz^x(\xi_0=x)=1$, and recursively for $n=1,2,\ldots$, 
\[
\Pozz^x[\xi_{J(n)} \in \cdot \mid \xi_{J(k)}, k < n]
  = \pi^{U_{[\xi_{J(n-1)}\cdot\e]}}(\cdot)
 \qquad {\rm if} \; \zeta'_n=1,
\]
and, for $\zeta'_n=0$,
\[
\Pozz^x[\xi_{J(n)} \in \cdot \mid \xi_{J(k)}, k<n] = 
\frac{ \Po^{\xi_{J(n-1)}}[\xi_{\eta_n} \in \cdot\;] - r_1
\pi^{U_{[\xi_{J(n-1)}\cdot\e]}}(\cdot)}{1-r_1}.
\]
Set $\kappa_0=0$, and define the times of success in the new 
Bernoulli process, 
\[
  \kappa_{m+1}=\min\{k > \kappa_m: \zeta'_k=1\}
 \]
 for $m\geq 1$. It is easy to see that, under 
$P^{\zeta'} \otimes  \Pozz^x$, the sequence $(\xi_{J(\kappa_m)},
m \geq 0)$ is a Markov chain, and $\xi_{J(\kappa_m)}$  
for $m\geq 1$ has ``piecewise stationary'' law,
i.e., of the form $\sum_{i \in \Z} a_i \pi^{U_i}$ with $a_i \geq 0,
 \sum_i a_i=1$. It follows that, starting $\xi_0$ from such a law,
the Markov chain is weakly lumpable and can be reduced to another
Markov chain 
on a smaller space, see~\cite{KemSn}, or~\cite{Jled} for a more 
modern account.  More precisely, we prove:
\begin{lmm}
\label{l_lumps}
Under $P^{\zeta'} \otimes  \Pozz^{U_0}$, the sequence 
$([\xi_{J(\kappa_m)}\cdot \e], m \geq 0)$ is 
a RWRE on $\Z$, with transition probabilities 
\[
Q_\omega(i,j)= P^{\zeta'} \otimes 
    \Pozz^{U_i}[\xi_{J(\kappa_1)} \in U_j].
\]
\end{lmm}

\noindent
\textit{Proof.} Under $P^{\zeta'} \otimes  \Pozz^x$, 
the transition density from~$x$ to~$y$ 
for the Markov chain $\xi_{J(\kappa_\cdot)}$ can be written 
as $a_j(x,\omega) \pi^{U_j}(dy)$,
with~$j$ such that $U_j \ni y$. Hence,
\begin{align*}
\lefteqn{P^{\zeta'} \otimes  \Pozz^{U_0}
\big[[\xi_{J(\kappa_1)}\cdot \e]=j_1, \ldots, 
 [\xi_{J(\kappa_m)}\cdot \e]=j_m\big]}
\\
&=
P^{\zeta'} \otimes  \Pozz^{U_0}\big[\xi_{J(\kappa_1)}\in U_{j_1}, 
\ldots, \xi_{J(\kappa_m)}\in U_{j_m}\big]\\
&=
\int\ldots \int \pi^{U_0}(dx_0) \prod_{k=1}^m a_{j_k}(x_{k-1},
\omega) \pi^{U_{j_k}}(dx_k)\\
&=
\prod_{k=1}^m\int  a_{j_k}(x_{k-1},\omega) 
\pi^{U_{j_{k-1}}}(dx_{k-1})  \times \int \pi^{U_{j_m}}(dx_m)\\
&=
\prod_{k=1}^m Q_\omega(j_{k-1},j_k) \times 1,
\end{align*} 
which ends the proof.
\qed

This result is the bridge between the two main processes considered
in this paper:
obviously, starting $\xi_0$ from the origin or distributed in the
interval~$U_0$ will not make any difference for the law of large
numbers.
But it is not quite enough to conclude the proof for the billiard.
%%%%%%%%%%%%%%%%%%%%%%%%%%%%%%%%%%
% Similarly to the construction in Section~\ref{s_proofs_rwre}, 
% we can
% formally define~$\Pozz$ in the following way. For all $i\in\Z$,
% consider a
% sequence $(Y_n^{(i)}, n=1,2,3,\ldots )$ of i.i.d.\ random points
% having the stationary distribution on~$U_i$. 
% Then, for all $n\geq 1$
% there exists a family of measures 
% (indexed by $x\in\partial\omega$)
% $\Po^{x,n}[(\eta,\xi,Y_n)\in\cdot\;]$ with the marginal
% $\Po^{x}[(\eta,\xi)\in\cdot\;]$ such that for all~$n$ and for
% all~$x\in U_i$
% \begin{equation}
% \label{coupling_Y}
%  \Po^{x,n}[\xi_{\eta_n}=Y_n^{(i)}]\geq r_1.
% \end{equation}
% Now, given $\xi_{J(n-1)}=y\in U_i$, the law of
% $(\eta_n, {\bar\xi}):=(\eta_n,\xi_{J(n-1)+1},\ldots,\xi_{J(n)})$
% under~$\Pozz^y$ is given by
% \begin{align*}
%  &\1{\zeta'_n=1}\Po^{y,n}[(\eta_n, {\bar\xi})\in\cdot \mid
% \xi_{\eta_n}=Y_n^{(i)}]\\
%  &\quad + \1{\zeta'_n=0}\Big( 
% \frac{\Po^{y,n}[\xi_{\eta_n}=Y_n^{(i)}]-r_1}{1-r_1}
% \Po^{y,n}[(\eta_n,
% {\bar\xi})\in\cdot \mid \xi_{\eta_n}=Y_n^{(i)}]\\
% &\qquad\qquad\qquad\qquad + 
% \frac{1-\Po^{y,n}[\xi_{\eta_n}=Y_n^{(i)}]}{1-r_1}
% \Po^{y,n}[(\eta_n,
% {\bar\xi})\in\cdot \mid \xi_{\eta_n}\neq Y_n^{(i)}]
% \Big)
% \end{align*}
% \textbf{this is wrong, since the conditional distribution of
% $\xi_{\eta_n}$ given $\zeta'_n=1$ need not be uniform (the
% conditional distribution of $Y$ is not necessarily so\dots)}
%%%%%%%%
%%%%%%%%%%%%
In the sequel we shall need the following two results about
hitting times of sets:
\begin{lmm}
\label{l_go_B}
For any $m\geq 0$ and arbitrary $B,F\subset\partial\omega$ we have
for $\IP$-almost all $\omega$
\begin{equation}
\label{eq_go_B}
 \Po^B[\text{there exists }k\leq m\text{ such that }\xi_k\in F]
   \leq m\frac{\pi(F)}{\pi(B)}.
\end{equation}
\end{lmm}

\noindent
\textit{Proof.}
Using reversibility, it is straightforward to obtain that
$\pi(B)\Po^B[\xi_k\in F]=\pi(F)\Po^F[\xi_k\in B]$, so
$\Po^B[\xi_k\in F]\leq\frac{\pi(F)}{\pi(B)}$. 
Using the union bound, we obtain~\eqref{eq_go_B}.
\qed

\begin{lmm}
\label{l_go_back_x}
There exist ${\tilde\gamma}_4,{\tilde\gamma}_5>0$ such that 
 for any $m\geq 0$, $H\geq 1$, and $x\in\partial\omega$ we have
for $\IP$-almost all $\omega$
\begin{equation}
\label{eq_go_back_x}
\Po^x[\text{there exists }k\leq m\text{ such that }(\xi_k-x)\cdot\e
< -H] \leq {\tilde\gamma}_4 m e^{-{\tilde\gamma}_5H^{1/2}}.
\end{equation}
\end{lmm}

\noindent
\textit{Proof.}
This fact would be a trivial consequence of Lemma~\ref{l_go_B} if
one starts from the stationary distribution on a set (of not too 
small measure) instead of starting from a single point. So, the
idea is the following: we first wait for the moment~$J(\kappa_1)$
(when the particle has the stationary distribution on~$U_j$ for
some random~$j$), and then note that it is likely that the particle
did not go too far to the left until this moment. Then, it is
already possible to apply Lemma~\ref{l_go_B}. 

Formally, we write
\begin{align*}
 \lefteqn{\Po^x[\text{there exists }k\leq m\text{ such that }
 (\xi_k-x)\cdot\e< -H]}\\
 &= E^{\zeta'}\Pozz^x[\text{there exists }
k\leq m\text{ such that }(\xi_k-x)\cdot\e < -H]\\
 &\leq P^{\zeta'}\Big[\kappa_1>\frac{H^{1/2}}{2N}\Big]\\
 & \quad
 +\Po^x\Big[\text{there exists }k\leq\frac{H^{1/2}}{2N}
\text{ such that } (\xi_k-\xi_{k-1})\cdot\e\leq -H^{1/2} \Big]\\
 &\quad + E^{\zeta'}\Pozz^x\Big[\kappa_1\leq\frac{H^{1/2}}{2N},
  (\xi_k-\xi_{k-1})\cdot\e >-H^{1/2}\text{ for all }
 k\leq \frac{H^{1/2}}{2N}, \\
 & \qquad\qquad\qquad\,\text{there exists }k\leq m
  \text{ such that
}(\xi_{k+J(\kappa_1)}-\xi_{J(\kappa_1)})\cdot\e<-H/2\Big].
\end{align*}
Now, the bound on the first term is straightforward (the random
variable~$\kappa_1$ has a geometric distribution with
parameter~$r_1$).
To estimate the second term, note that
 for any~$h>0$ one has
\begin{equation}
\label{jump_neg}
  \Po^x[(\xi_1-x)\cdot\e < -h] 
= \intl_{(y-x)\cdot \e < -h} e^{-\lambda (y-x)\cdot \e} K(x,y) dy
\leq e^{-\lambda h}.
\end{equation}
To deal with the third term, recall that the law of $\{\xi_j, j >
J(\kappa_1)\}$ conditional on~$\xi_{J(\kappa_1)}$ does not depend
%%on~$\mathcal{F}_{J(\kappa_1)}$. 
on $(\xi_m; m \leq J(\kappa_1))$.
Then,
on the event 
\[
 \Big\{\kappa_1\leq\frac{H^{1/2}}{2N},
  (\xi_k-\xi_{k-1})\cdot\e >-H^{1/2}\text{ for all }
 k\leq \frac{H^{1/2}}{2N}\Big\}
\]
we have $[\xi_{J(\kappa_1)}\cdot\e]\geq -H/2$. So, one can estimate
the third term using Lemma~\ref{l_go_B}, and conclude the proof of
Lemma~\ref{l_go_back_x}.
\qed

\medskip

Now, to prove Theorem~\ref{t_kn_lln}, the idea is to construct an
``induced'' RWRE, then apply the results of
Section~\ref{s_results_rwre}, and then recover the LLN for the
original billiard. To apply this approach, we need a few estimates
on displacement probabilities 
that we derive in the following lines.
Consider some (large) integer~$L$ (to be chosen later) and observe
that, since $\kappa_n$ is a sum of~$n$ i.i.d.\ geometric random
variables, we can find some large $r_1$ such that, for all $n$,
\begin{equation}
\label{NL4} 
 P^{\zeta'}[\kappa_{n}\leq 2r_1^{-1}n]\geq 1-C_7e^{-C_8n}.
\end{equation}
 So, by Lemma~\ref{l_go_B}, \eqref{NL4}, and using also the fact
that $J(\kappa_{L^4})\leq N\kappa_{L^4}$, we obtain for
$\IP$-almost all $\omega$ that for every $m\geq 1$ it holds that
\begin{align}
 \lefteqn{E^{\zeta'}\Pozz^{U_0}\big[[\xi_{J(\kappa_{L^4})}\cdot\e]
   \in [-(m+1)L,-mL]\big]}~~~~~~\nonumber\\
&\leq E^{\zeta'}\Pozz^{U_0}\big[\text{there exists }
 k\leq 2mNr_1^{-1}L^4\nonumber\\
&\qquad\qquad\qquad\qquad
\text{ such that }\xi_k\cdot\e\in [-(m+1)L,-mL]\big]
\nonumber\\
&\quad+ P^{\zeta'}[\kappa_{L^4}>2mr_1^{-1}L^4]\nonumber\\
 &\leq C_9mL^4e^{-C_{10}mL} +C_7e^{-C_8mL^4}.
\label{displace_1}
\end{align}
Also, we use Lemma~\ref{l_go_out} (applied to ${\tilde
F}^\omega(-L,2L)$), \eqref{NL4}, Lemma~\ref{l_go_back_x}, and the
fact that $J(\kappa_{L^4})\geq L^4$, to
obtain that for $\IP$-almost all~$\omega$
\begin{align}
 \lefteqn{E^{\zeta'}\Pozz^{U_0}\big[[\xi_{J(\kappa_{L^4})}\cdot\e]
\geq L\big]}\nonumber\\ 
&\geq E^{\zeta'}\Pozz^{U_0}[\kappa_{L^4}\leq
2r_1^{-1}L^4,\tau(-L,2L)<L^4, \xi_k\cdot\e>-L,
\nonumber\\ 
&\qquad\qquad\qquad
(\xi_{k+\tau(-L,2L)}-\xi_{\tau(-L,2L)})\cdot\e>-L
\text{ for all }k\leq 2Nr_1^{-1}L^4]\nonumber\\ 
&\geq
1-C_7e^{-C_8L^4}-
{\tilde\gamma}_1e^{-{\tilde\gamma}'_1L/27}
-2Nr_1^{-1}L^4e^{-\lambda L}
-2{\tilde\gamma}_4Nr_1^{-1}L^4e^{-{\tilde\gamma}_5L^{1/2}}.
\label{displace_2}
\end{align}
So, from~\eqref{displace_1} and~\eqref{displace_2} we obtain
\begin{align}
 \lefteqn{E^{\zeta'}\Eozz^{U_0}[\xi_{J(\kappa_{L^4})}\cdot \e]}
\nonumber\\
  &\geq L(1-C_7e^{-C_8L^4}-
{\tilde\gamma}_1e^{-{\tilde\gamma}'_1L/27}
-2Nr_1^{-1}L^4e^{-\lambda L}
-2{\tilde\gamma}_4Nr_1^{-1}L^4e^{-{\tilde\gamma}_5L^{1/2}})
 \nonumber\\ 
 &\quad -L(C_7e^{-C_8L^4}+
{\tilde\gamma}_1e^{-{\tilde\gamma}'_1L/27}
+2Nr_1^{-1}L^4e^{-\lambda L}
+2{\tilde\gamma}_4Nr_1^{-1}L^4e^{-{\tilde\gamma}_5L^{1/2}})
 \nonumber\\ 
 &\quad -\sum_{m=1}^\infty mL(C_7e^{-C_8mL^4}+C_9mL^4e^{-C_{10}mL})
\nonumber\\ 
&>1 \label{drift_billiard}
\end{align}
if~$L$ is large enough. 

\medskip

\noindent
\textit{Proof of Theorem~\ref{t_kn_lln}.} 
Under the law $E^{\zeta'}\Pozz^{U_0}(\cdot)$,
the process~$S$ defined by
\[
S_n:=[\xi_{J(\kappa_{L^4n})}\cdot\e], \qquad n \geq 0,
\]
is a RWRE on~$\Z$ in a (stationary ergodic) 
environment given by the tube~$\omega$. 
Let us choose~$L$ such 
that~\eqref{drift_billiard} holds.
In this case, \eqref{displace_1} and~\eqref{drift_billiard} show
that the process~$S$ has
uniformly positive drift, and its jumps to the left have uniformly
exponential tail. 

Now, to apply the results of Section~\ref{s_results_rwre} to 
 the process~$S$, we need
to check that it verifies Conditions~E, C, D. First, it
is straightforward to obtain that Condition~E holds. To check
Condition~C, first recall that (cf.\ e.g.\ formula~(54)
of~\cite{CPSV1}) 
that there exists~${\tilde\gamma}_6>0$
(depending only on~$\widehat M$) such that for
all~$x\in\partial\omega$ and all~$h\geq 1$
\begin{equation}
\label{long_jump}
 \Po^x\big[|(\xi_1-x)\cdot\e|>h\big] \leq
 {\tilde\gamma}_6 h^{-(d-1)}.
\end{equation}
Now, observe that
%, for~$\zeta'$ such that $\kappa_1=j$, we have
\begin{equation}
\label{=j}
 [\xi_{J(j)}\cdot\e] = [\xi_{J(j-1)}\cdot\e]
\text{ on } \{j=\kappa_n \text{ for some }n\}
\end{equation}
and, for~$i$ such that $i\neq\kappa_n$ for all~$n$,
\begin{align}
 \lefteqn{E^{\zeta'} \big(\Pozz^{U_0}\big[|[\xi_{J(i)}\cdot\e]-
  [\xi_{J(i-1)}\cdot\e]| \geq h \big]  \mid \kappa_{L^4}=j \big)}
\nonumber\\
 &= E^{\zeta'} \big(\Pozz^{U_0}\big[|[\xi_{J(i)}\cdot\e]-
  [\xi_{J(i-1)}\cdot\e]| \geq h \big]  \mid \zeta'_i=0 \big)
\nonumber\\
&\leq \frac{1}{P^{\zeta'}[\zeta'_i=0]}
  \Po^{U_0}\big[|[\xi_{J(i)}\cdot\e]-
  [\xi_{J(i-1)}\cdot\e]| \geq h \big]
\nonumber\\
&\leq C_{11}h^{-2}
\label{<j}
\end{align}
since $d\geq 3$, recall~\eqref{long_jump}. 
Then, write using~\eqref{=j} and~\eqref{<j}
\begin{align}
\lefteqn{E^{\zeta'} \Pozz^{U_0}\big[|[\xi_{J(\kappa_{L^4})} \cdot
\e]|\geq s\big]}
\nonumber\\
 &= \sum_{j=1}^\infty P^{\zeta'}[\kappa_{L^4}=j]
  E^{\zeta'}\big(\Pozz^{U_0}\big[|[\xi_{J(\kappa_{L^4})}\cdot\e]|
\geq s\big]\mid \kappa_{L^4}=j\big)
\nonumber\\
&\leq \sum_{j=1}^\infty P^{\zeta'}[\kappa_{L^4}=j]
  E^{\zeta'}\big(\Pozz^{U_0}\big[\text{there exists }i\leq j
\text{ such that }
\nonumber\\
& \qquad\qquad\qquad\qquad\qquad\qquad\qquad
|[\xi_{J(i)}\cdot\e]-
  [\xi_{J(i-1)}\cdot\e]|
\geq s/j\big]\mid \kappa_{L^4}=j\big)
\nonumber\\
&\leq \sum_{j=1}^\infty P^{\zeta'}[\kappa_{L^4}=j]
jC_{11}\Big(\frac{s}{j}\Big)^{-2}
\nonumber\\
&= C_{11} s^{-2} \sum_{j=1}^\infty j^3P^{\zeta'}[\kappa_{L^4}=j]
\nonumber\\
&= C_{12} s^{-2}.
\label{oc_kappa_L4}
\end{align}
Abbreviating $P^*[\cdot]:=E^{\zeta'} \Pozz^{U_0}[\cdot]$, we see
that~\eqref{oc_kappa_L4} is equivalent to $P^*[|S_1|\geq s]\leq
C_{12} s^{-2}$.
This means that Condition~C holds for the process~$(S_n,n\geq 0)$.
% To see that Condition~C holds for the process~$(S_n,n\geq 0)$, it
% is
% enough to observe that~\eqref{oc_kappa_1} is satisfied also for 
% $|[\xi_{J(\kappa_i)}\cdot\e]-[\xi_{J(\kappa_{i-1})}\cdot\e]|$
% on the place of $|[\xi_{J(\kappa_1)}\cdot\e]|$, and use the
% inclusion 
% \begin{align*}
%  \{|[\xi_{J(\kappa_{L^4})} \cdot\e]|\geq s\}
% &\subset
%  \{\text{there exists }i\leq L^4\\
%  &\qquad\qquad\text{ such that }
%  |[\xi_{J(\kappa_i)}\cdot\e]
%    -[\xi_{J(\kappa_{i-1})}\cdot\e]|\geq s/L^4\}.
% \end{align*}

Now, we show that Condition~D holds 
for the process~$(S_n,n\geq 0)$.
First, using~\eqref{displace_1} and~\eqref{drift_billiard} and
Condition~C, for large enough~$\vr_0$
one obtains by a straightforward computation that there exist 
small enough ${\tilde\gamma}_7, {\tilde\gamma}_8>0$ such that 
for all $y\in\Z$ and all~$\vr\geq\vr_0$
\[
 E^*(e^{-{\tilde\gamma}_7 S^\vr_{n+1}+{\tilde\gamma}_8(n+1)}
      -e^{-{\tilde\gamma}_7S^\vr_n+{\tilde\gamma}_8n}
   \mid S^\vr_n=y) \leq 0,
\]
so that $(e^{-{\tilde\gamma}_7S^\vr_n+{\tilde\gamma}_8n},n\geq 0)$
 is a positive
supermartingale. So, for some positive constants $C_{13}, C_{14}$
it holds that for all~$k\geq 1$
\[
 P^*[\text{there exists }n\text{ such that }S^\vr_n\leq -k] <
   C_{13}e^{-{\tilde\gamma}_7k}
\]
and
\[
 E^*\sum_{n=0}^\infty \1{S^\vr_n=0} < C_{14}.
\]
So, Condition~D holds with
$g_1(k)=C_{13}C_{14}e^{-{\tilde\gamma}_7k}$.
 This means that we can use Theorem~\ref{t_rw_lln}
for the process~$S$.

Now, 
it remains to deduce the LLN for the random billiard with drift
from the LLN for the random walk in random environment.
It is done by a standard argument that we sketch in the following
lines. First, observe that, just in the same way
as~\eqref{oc_kappa_L4} one can prove a slightly more general fact:
for any~$n$
\begin{equation}
\label{gener_kappa_L4}
 E^{\zeta'} \Pozz^{U_0}\Big[\max_{m\in
\big[J(\kappa_{L^4n}),J(\kappa_{L^4(n+1)})\big)}
|[\xi_m-\xi_{J(\kappa_{L^4n})} \cdot
\e]|\geq s\Big] \leq C_{15} s^{-2}.
\end{equation}
Then, since the limit of $n^{-1}S_n$ exists and is finite, there
exists $u\in(0,\infty)$ such that
\[
 u=\lim_{n\to\infty} 
 \frac{\xi_{J(\kappa_{L^4n})}\cdot\e}{J(\kappa_{L^4n})}.
\]
We then use~\eqref{gener_kappa_L4} to obtain that,
for $m\in
\big[J(\kappa_{L^4n}),J(\kappa_{L^4(n+1)})\big)$
\begin{align*}
 \frac{1}{n}|\xi_{J(\kappa_{L^4n})} - \xi_m| &\leq
\frac{1}{n}
 \max_{m\in
\big[J(\kappa_{L^4n}),J(\kappa_{L^4(n+1)})\big)}
|[\xi_m-\xi_{J(\kappa_{L^4n})} \cdot
\e]|\\
&\to 0 \quad \text{ a.s., as }n\to\infty,
\end{align*}
and this permits us to conclude the proof of 
Theorem~\ref{t_kn_lln}.
\qed

\end{document}